\documentclass[11pt]{amsart}

\usepackage{quiver}
 
\usepackage{amssymb,graphicx, amsmath, amsthm,MnSymbol, array}
\usepackage{mathtools}
\usepackage{extarrows}

\graphicspath{ {./images/} }
\usepackage{tikz-cd}
\definecolor{darkgreen}{rgb}{0,0.4,0}
\DeclareMathAlphabet\mathbfcal{OMS}{cmsy}{b}{n}

\usepackage{todonotes, comment}

\usepackage{floatflt}

\usepackage{hyperref}

\usepackage{ulem}
\usepackage{mathbbol}

\usepackage[inline]{enumitem}

\usepackage[mathscr]{eucal}

\usepackage[font=small]{caption}

\usepackage{etoolbox}

\makeatletter
\def\l@subsection{\@tocline{2}{0pt}{2pc}{5pc}{}}
\makeatother

\def\m{{\mathfrak m}} %%% m

\def\k{\mathbb R} %%% R
\def\k{{\kappa}}

\def\O{{\mathcal{O}}}

\def\PP{f}

\def\PP{\overline{\PP}}

\begin{document}

\newtheorem{theorem}{Theorem}[section]
\newtheorem{maintheorem}[theorem]{Main Theorem}
\newtheorem{lemma}[theorem]{Lemma}
\newtheorem{proposition}[theorem]{Proposition}
\newtheorem{corollary}[theorem]{Corollary}
\newtheorem{problem}[theorem]{Problem}
\newtheorem{construction}[theorem]{Construction}

\theoremstyle{definition}
\newtheorem{defi}[theorem]{Definitions}
\newtheorem{definition}[theorem]{Definition}
\newtheorem{notation}[theorem]{Notation}
\newtheorem{remark}[theorem]{Remark}
\newtheorem{example}[theorem]{Example}
\newtheorem{question}[theorem]{Question}
\newtheorem{Comment}[theorem]{Comment}
\newtheorem{comments}[theorem]{Comments}

\renewcommand{\thedefi}{}

\long\def\alert#1{\smallskip{\hskip\parindent\vrule%
\vbox{\advance\hsize-2\parindent\hrule\smallskip\parindent.4\parindent%
\narrower\noindent#1\smallskip\hrule}\vrule\hfill}\smallskip}

\def\ann{\mathrm{ann}}
\def\red{{{\rm red}}}
\def\ff{\mathfrak}
\def\tf{torsion-free}
\def\Spec{\mbox{\rm Spec}}
\def\Proj{\mbox{\rm Proj }}
\def\hgt{\mbox{\rm ht }}
\def\type{\mbox{ type}}
\def\Hom{\mbox{ Hom}}
\def\rank{\mbox{rank}}
\def\Ext{\mbox{ Ext}}
\def\Tor{\mbox{ Tor}}
\def\ker{{{\rm{ker \: } }}}
\def\Max{\mbox{\rm Max}}
\def\Min{\mbox{\rm Min}}
\def\End{\mbox{\rm End}}
\def\xpd{\mbox{\rm xpd}}
\def\Ass{\mbox{\rm Ass}}
\def\emdim{\mbox{\rm emdim}}
\def\epd{\mbox{\rm epd}}
\def\repd{\mbox{\rm rpd}}
\def\ord{\mbox{\rm ord}}
\def\Clop{\mbox{\rm Clop}}
\def\fod{\mbox{\rm fod}}
\def\od{\mbox{\rm od}}

\def\Id{\mbox{\rm Id}}

\def\Bis{\mbox{\rm Bis}}

\def\htt{\mbox{\rm ht}}

\def\DD{{\mathcal D}}
\def\EE{{\mathcal E}}
\def\FF{{\mathcal F}}
\def\GG{{\mathcal G}}
\def\HH{{\mathcal H}}
\def\II{{\mathcal I}}
\def\LL{{\mathcal L}}
\def\MM{{\mathcal M}}
\def\PP{{\mathcal P}}

\def\k{\mathbb{k}}

\def\cal{\mathcal}

\def\BPR{{\bf B^P_\textit{R}}}
\def\PBR{{\bf P^B_\textit{R}}}
\def\ABR{{\bf A^B_\textit{R}}}
\def\BAR{{\bf B^A_\textit{R}}}
\def\PAR{{\bf P^A_\textit{R}}}
\def\APR{{\bf A^P_\textit{R}}}
\def\ASR{{\bf A^S_\textit{R}}}
\def\PSR{{\bf P^S_\textit{R}}}
\def\SAR{{\bf S^A_\textit{R}}}
\def\SPR{{\bf S^P_\textit{R}}}

\title{Geometry of Spaces of Subrings}

%    author one information
\author{Laura Cossu}
\address{\parbox{\linewidth}{Dipartimento di Matematica e Informatica\\ Università degli Studi di Cagliari\\ Palazzo delle Scienze, Via Ospedale 72, 09124 Cagliari, Italy}}
\email{laura.cossu3@unica.it}

%    author two information
\author{Bruce Olberding}
\address{ \parbox{\linewidth}{Department of Mathematical Sciences\\
New Mexico State University\\ 1290 Frenger Mall, Las Cruces, New Mexico 88003-8001, USA}}
\email{bruce@nmsu.edu}

%    author three information
\author{Giovanni Secreti}
\address{\parbox{\linewidth}{Department of Mathematical Sciences\\
New Mexico State University\\ 1290 Frenger Mall, Las Cruces, New Mexico 88003-8001, USA}}
\email{gsecreti@nmsu.edu}

\begin{abstract} Motivated by examples such as the local rings of projective integral schemes and the Zariski-Riemann space of a field, we develop a geometry for continuously indexed families of subrings
of a commutative ring $R$. The indexing spaces are Stone spaces, and continuity is taken with respect to the Zariski topology on the space of subrings. To each such family we associate a subsheaf of
the constant sheaf with value $R$ whose stalks are the indexed subrings. This gives a \textit{patch ringed space}. Its global sections are locally constant $R$-valued functions whose value at each
point belongs to the corresponding subring. The ring of global sections, with this representation, is the associated \textit{patch algebra}.

Our main result establishes an algebra-geometry correspondence between these continuously indexed families of subrings and patch ringed spaces that can be expressed as a category equivalence, with each of the two categories dually equivalent to the category of patch algebras. Thus the indexing space, its sheaf, and the indexed subrings can be recovered functorially from the associated patch algebra, which serves as a coordinate ring for this geometry. We describe kernels and factorizations of patch algebra morphisms in geometric terms and characterize their
injectivity and surjectivity. We include applications to Pierce duality and Bhatt-Scholze $w$-localizations of  integral separated quasi-compact schemes. 
\end{abstract}


\maketitle

\begingroup
\small
\tableofcontents
\endgroup

\section{Introduction}

Throughout this article, all rings are commutative, and we assume $R$ is a ring with $0 \ne 1$.
Spaces of subrings arise in many contexts, such as the Zariski-Riemann space of valuation rings of a field, the collection of local rings of a projective integral scheme, or the space of integrally closed subrings of a given ring. These and related spaces have been the subject of much recent study; see, for example, {\cite{CamFin26, Fin, FFL, FFS, FinSpir20, HLOT24, SpiritoLoc, SpiritoIso22}}. 
Motivated by these applications, our purpose is to formulate a geometry for spaces of subrings of a given ring. Here, geometry means viewing such a space not only as a topological space but also as a space equipped with a structure sheaf, a ``coordinate ring,'' sections, and local data that {\it functorially} encode its essential information.

In \cite{CosOlb}, the first two authors proved a realization theorem for spaces of subrings of a commutative ring. Given  a ring $R$, suppose that $X$ is a spectral subspace of the space of subrings of $R$, endowed with the Zariski topology. Then there exists a commutative ring $A$, called a {\it patch algebra}, that efficiently encodes the rings in $X$ as factor rings of $A$. When $R$ is an integral domain, this efficiency can be described as follows: the minimal spectrum of $A$ is homeomorphic to $X$  with its patch topology, and the rings in $X$ are the factor rings of $A$ by its minimal prime ideals. As discussed in \cite{CosOlb}, these ideas can be applied to construct patch algebras associated with the spaces mentioned above, as well as with spaces of local subrings of a ring, among many other examples.

This construction works more generally. Let $\Sigma(R)$ denote the space of subrings of $R$ endowed with the Zariski topology. If $Y$ is a Stone space (i.e., a zero-dimensional compact Hausdorff space, also known as a Boolean or profinite space) and $f\colon\; Y\rightarrow\Sigma(R)$ is a continuous map, then there is a patch algebra $A(f)$ encoding the rings in $f(Y)$ as factor rings of $A(f)$. In \cite{CosOlb}, such a continuous map $f$ is called a {\it patch bundle}. A patch bundle may be viewed as a continuously indexed family of subrings of $R$, with indexing space $Y$. Thus, the space of subrings $f(Y)$ need not itself be spectral for the construction of its associated patch algebra; it only needs to admit a continuous indexing by a Stone space.

The way in which a patch algebra $A(f)$ captures the structure of a patch bundle $f$ of subrings suggests that a patch algebra should be regarded as a coordinate ring for a space of subrings. We develop this analogy in Section~\ref{sec: patch algebras and patch bundles} by defining a category of patch algebras for $R$ whose objects include the structural data needed to recover their associated patch bundles. We show that these objects are precisely the patch algebras arising from patch bundles. In keeping with the coordinate ring analogy, we show that there is a duality between the category of patch bundles over $R$ and the category of patch algebras for $R$, so that patch algebras are indeed the algebraic counterparts of patch bundles. Furthermore, morphisms of patch algebras behave like those of coordinate rings: restriction to a closed subspace of the indexing Stone space induces a surjection of patch algebras, and every morphism of patch algebras factors as a surjection followed by an injection (see Section~\ref{sec: morphisms and kernels}).

 A patch bundle of subrings of $R$ contains more information than its underlying topology, since the indexed subrings themselves are essential to our constructions. In \cite{CosOlb}, this additional structure is encoded by a presheaf of subrings of $R$ on a Stone space. Such presheaves are called {\it patch presheaves for~$R$}, and their category is equivalent to the category of patch bundles over $R$. 
 This gives a different model of patch bundles and patch algebras and suggests a geometric interpretation involving presheaves, sections, and local data.

 In the present article, we provide such an interpretation by introducing patch ringed spaces, which are ringed Stone spaces $(Y,\mathcal O_Y)$ whose structure sheaf $\O_Y$ is a subsheaf of the constant sheaf~$\underline R_Y$ on $Y$ (see Section~\ref{sec: patch ringed spaces}). Our main result is that the relationship between these different structures is functorial and expressible as a triangle of category equivalences.  
 
\medskip

{\noindent}{\bf Main Theorem} (Theorems~\ref{thm: bundle ringed equivalence} and \ref{ASR}){\bf .} \textit{For a fixed ring $R$, the following diagram of categories commutes up to natural isomorphism, and all its arrows are equivalences:}
\[\begin{tikzcd}
	\mbox{Patch bundles} &&&& \mbox{Patch ringed spaces} \\
	\\
	&& {\mbox{(Patch algebras)}^{\rm op}}
	\arrow["{\mathcal S}", from=1-1, to=1-5]
	\arrow["{\mathcal A}"', from=1-1, to=3-3]
	\arrow["\Gamma", from=1-5, to=3-3]
\end{tikzcd}.\]
\textit{Here, $\Gamma$ is the global-sections functor. For a patch bundle $f\colon\; Y\to \Sigma(R)$ and its associated patch ringed space $\mathcal S(f)=(Y,\mathcal O_f)$, there are natural isomorphisms  
\[\Gamma(Y,\mathcal O_f)\cong A(f)
 \: {\mbox{ and }} \: 
(\mathcal O_f)_y\cong f(y),\]
 with the stalks of $\O_f$ identified with the subrings of $R$ in the image of $f$.}

\medskip
 
In this way, we obtain a geometry of subrings that is an analogue of the geometry of affine  schemes. We summarize this analogy in the following table, where  $\mathfrak m_y=\{U\in\Clop(Y):y\notin U\}$ is the maximal ideal of $\Clop(Y)$ corresponding to $y$ via Stone duality.

\begin{center}
\renewcommand{\arraystretch}{1.4}
 \begin{tabular}{|c|c|}
\hline
\textbf{Affine geometry} & \textbf{Patch geometry of subrings} \\
\hline
affine scheme $X$ 
& patch bundle $f\colon\; Y\to \Sigma(R)$ \\
\hline
structure sheaf $\mathcal O_X$
& subsheaf $\mathcal O_Y \subseteq \underline{R}_Y$ \\
\hline
coordinate ring $A=\Gamma(X,\mathcal O_X)$
& patch algebra $A(f)=\Gamma(Y,\mathcal O_Y)$ \\
\hline
point $x\in X$
& point $y\in Y$ \\
\hline
stalk $=$ local ring $\mathcal O_{X,x}$
& stalk $=$ subring $f(y)\subseteq R$ \\

$=$ localization of $A$
& $=$ factor ring $A(f)/\mathfrak m_yA(f)$ \\
\hline
closed subscheme $\operatorname{Spec}(A/I)\hookrightarrow X$
& closed subspace $Y'\subseteq Y$ \\
\hline
surjection $A\to A/I$
& surjection $A(f)\to A(f|_{Y'})$ \\
\hline
\end{tabular}
\end{center}
\smallskip

\smallskip

We apply these ideas to several running examples through the article, including the trivial, but instructive case of finitely many subrings; a chain of overrings of the ring of integers; the Zariski-Riemann space of a field; and integral separated quasi-compact schemes.  We also show in two appendices how these constructions connect to Pierce duality for commutative rings (Theorem~\ref{Pierce duality restriction}) and the work of Bhatt-Scholze on $w$-localizations of integral separated  quasi-compact schemes (Theorem~\ref{thm: scheme w-localization}).

\subsection*{Plan of the paper}
In Section~\ref{sec: preliminaries} we collect useful preliminaries on the Zariski and patch topologies on $\Sigma(R)$ and on Specker algebras. Section~\ref{sec: bundles and presheaves} recalls from \cite{CosOlb} the notions of patch bundles and patch presheaves, together with their categorical equivalence, and develops several examples that will be used throughout the article. We define the category of patch algebras for $R$ in
Section~\ref{sec: patch algebras and patch bundles}, where we also prove its duality with the category of patch bundles. In Section~\ref{sec: morphisms and kernels} we study morphisms of patch algebras and their kernels, and characterize injectivity and surjectivity in terms of the associated patch bundles. Section~\ref{sec: patch ringed spaces} introduces patch ringed spaces and establishes
their equivalence with patch bundles and patch algebras through sheafification and global sections. Appendices~\ref{sec: Pierce} and \ref{sec: w localization} describe the connections with Pierce
duality and $w$-localization, respectively.

\section[Preliminaries]{Preliminaries}\label{sec: preliminaries}
We recall some fundamental definitions and results from \cite{CosOlb} that will serve as a foundation for this work. All rings and algebras considered in this paper are commutative with identity. All subrings share the same identity, and all ring and algebra homomorphisms preserve identity. 
Throughout the rest of the article, we assume $R$ is a ring with $0 \ne 1$ and $\Sigma(R)$ is the set of subrings of $R$. 
 We allow the empty Stone space, the trivial Boolean algebra with $0=1$, and the zero patch algebra.

Stone duality is the dual equivalence between the category of Boolean algebras and Boolean homomorphisms and the category of Stone spaces and continuous maps, given on objects by $B\mapsto\Max(B)$ and $Y\mapsto\Clop(Y)$; see \cite[Corollary~II.4.4, p.~71]{Joh86}. A Boolean homomorphism $h\colon\; B_1\to B_2$ induces the continuous map $\Max(B_2)\to\Max(B_1)$ defined by $\mathfrak m\mapsto h^{-1}(\mathfrak m)$, and a continuous map of Stone spaces induces a Boolean homomorphism in the opposite direction by inverse image on clopen subsets. 

For a Boolean algebra $B$ and $e\in B$, write
$$U_e=\{\mathfrak m\in\Max(B):e\notin\mathfrak m\}.$$
We use this notation throughout the paper. Under Stone duality, the assignment $e\mapsto U_e$ identifies $B$ with $\Clop(\Max(B))$. When $B=\Clop(Y)$ for a Stone space $Y$, we also write $U_e$ for the corresponding clopen subset of $Y$ under the canonical homeomorphism $Y\cong\Max(B)$. 

\subsection[Spaces of subrings]{Spaces of subrings}\label{subsec: spaces of subrings} Let $\Sigma(R)$ denote the set of all subrings of $R$. The Zariski topology on $\Sigma(R)$ has basic open sets
$$
\mathcal U(r_1,\ldots,r_n)=\{T\in\Sigma(R):r_1,\ldots,r_n\in T\}.
$$
With this topology $\Sigma(R)$ is a spectral space.\footnote{The Zariski topology on $\Sigma(R)$ coincides with the Scott topology on the algebraic lattice of subrings of $R$, placing this space within the general framework of spectral spaces of subobjects recently studied by Campanini and Finocchiaro \cite[Proposition~4.6 and Corollary~5.1]{CamFin26}.} Besides the Zariski topology, we can endow $\Sigma(R)$ with the {\it patch topology}, a Hausdorff refinement of the Zariski topology. A basis for the patch topology of $\Sigma(R)$ consists of the sets
$$
\mathcal U(r_1,\ldots,r_n)\setminus \bigcup_{j=1}^{m}\mathcal U(s_j), \:\: n,m\geq 0,$$ where $r_i,s_j\in R$, with the convention that $\mathcal U(r_1,\ldots,r_n)=\Sigma(R)$ when $n=0$. In the patch topology, $\Sigma(R)$ is a Stone space, and we will denote it as ${\Sigma}(R)^{\rm patch}$.

We recall from \cite[Definition 3.2 and Proposition 3.1]{CosOlb} that a {\it patch space} for $R$ is a closed subspace of the space ${\Sigma}(R)^{\rm patch}$, or, equivalently, a spectral subspace of $\Sigma(R)$. As a closed subspace of a Stone space, every patch space is Stone.

If $R$ is a ring and $S\subseteq R$ a subring, examples of patch spaces include the space of all subrings of $R$ containing $S$, the space of all rings between $S$ and $R$ integrally closed in $R$, and, when $R$ is a field, the Zariski-Riemann space consisting of the valuation rings of $R$ that contain $S$. We refer to \cite[Example 3.3]{CosOlb} for further examples. 

\subsection[Specker algebras]{Specker algebras}\label{subsec: Sepecker algebras}
In Section~\ref{sec: patch algebras and patch bundles} we will define patch algebras as rings equipped with specified embeddings into what are called Specker algebras in \cite[Definition~2.3]{BMMO}. To define these algebras, we recall from \cite[p.~185]{BMMO} the notion of a faithful generating algebra of idempotents of an $R$-algebra $S$, where $R$ is a ring. 
A nonzero idempotent $e \in S$ is {\it faithful} if $re=0$ implies $r=0$ for every $r\in R$. Let $B$ be a Boolean subalgebra of the Boolean algebra of idempotents of $S$, and suppose that $B$ generates $S$ as an $R$-algebra. Then $ B$ is a {\it faithful
generating algebra of idempotents of $S$} if each nonzero $e \in B$ is faithful. 

A {\it Specker $R$-algebra} is an $R$-algebra that has a faithful generating algebra $B$ of idempotents. Thus  \[S=\{t_1e_1 + \cdots + t_n e_n : n>0, e_i \in B\text{ and }t_i \in R \text{ for every }i\in \{1,\ldots, n\}\}.\]
Although the choice of $B$ need not be unique, any two such Boolean algebras are isomorphic \cite[Example~3.4 and Theorem~3.5]{BMMO}.
It follows that $B$ and $R$ uniquely determine $S$ up to isomorphism, in the sense that every Specker $R$-algebra is isomorphic as an $R$-algebra to the ring $R[B]$ whose elements are \textit{formal} sums of the form $\sum_i r_i{e_i}$, with each $r_i \in R$ and $e_i \in B$; see \cite[Definition~2.4 and Theorem~2.7]{BMMO} for details on this construction. Here the top element $1_B$ of $B$ is $1_S$, which is the image of $1_R$ under the structure map $R\to S$.

Every element $a\in R[B]$ admits a unique full orthogonal decomposition with distinct coefficients \cite[Lemma 2.1]{BMMO}:
    \[ a=\sum_{i=1}^n t_i e_i, \text{ where the }t_i\in R\text{ are distinct and }e_i\in B.\]
Here, by full orthogonal decomposition we mean that the $e_i$ in the above equality are all nonzero, pairwise orthogonal, and that $\bigvee_{i=1}^n e_i=1$. For the zero algebra, we regard the empty sum as the full orthogonal decomposition of zero. We use $\fod(a)$ to denote this decomposition and, if $a\ne 0$, $\od(a)$ to denote the orthogonal decomposition with nonzero distinct coefficients obtained from $\fod(a)$ by removing (if it exists) the summand with zero coefficient.

\section[Patch bundles and their presheaves]{Patch bundles and their presheaves}\label{sec: bundles and presheaves}

 We recall from \cite{CosOlb} the categories of patch bundles and patch presheaves, together with their categorical equivalence. We also introduce new examples that enrich preexisting the theory and will serve as running examples guiding the discussion throughout the paper.

\subsection[Patch bundles]{Patch bundles}\label{subsec: patch bundles} A {\it patch bundle over the ring $R$} is a continuous map $f\colon\;  Y \rightarrow \Sigma(R)$ from a Stone space~$Y$ to $\Sigma(R)$.
With a slight abuse of notation, we will sometimes identify the patch bundle $f\colon\; Y \rightarrow \Sigma(R)$ with the surjective continuous map $Y\to X$, where $X=f(Y)$. 

 As in \cite[Definition 3.11]{CosOlb}, a {\it morphism of patch bundles over $R$}
\[(f_1\colon\; Y_1\to {\Sigma(R)})\to (f_2\colon\; Y_2\to {\Sigma(R)})\]
is a continuous map $g\colon\; Y_1\to Y_2$ such that
\[f_2(g(y_1))\subseteq f_1(y_1)\text{ for every }y_1\in {Y_1}.\] We denote by  $\mathsf{PBun}_R$ the category of patch bundles over $R$ with their morphisms.

Let $X\subseteq\Sigma(R)$ be a subset whose Zariski subspace topology is spectral. The inclusion $\iota_X\colon\; X^{\rm patch}\to\Sigma(R)$ is a patch bundle, where $X^{\rm patch}$ carries its intrinsic patch topology. We call this the {\it tautological patch bundle} associated with $X$. When $X$ is a patch space, its intrinsic patch topology agrees with the topology induced by $\Sigma(R)^{\rm patch}$.

\medskip

We give four examples of patch bundles, each of which we will revisit in later sections. 
Before proceeding, let us recall that a continuous map between spectral spaces is {\it spectral} if preimages of quasi-compact open sets are quasi-compact. Thus, a patch bundle $f\colon\; Y\to \Sigma(R)$ is spectral if and only if, for each $r\in R$, the set $\{y\in Y : r\in f(y)\}$ is quasi-compact, or equivalently, clopen.

\begin{example} \label{finite 1} Let $X$ be a finite collection of subrings of the ring $R$. Then the patch topology on $X$ is the discrete topology, 
and the tautological patch bundle $X^{\rm patch} \rightarrow X$ is a spectral map. However, patch bundles with finite image need not be spectral. To see this, let $X=\{S,T\}$,
where $S\subsetneq T\subseteq R$. For any Stone space
$Y$ and any open subset $U\subseteq Y$, define
$f:Y\to X$ by
\[
f(y)=
\begin{cases}
T & \text{if } y\in U,\\
S & \text{if } y\notin U.
\end{cases}
\] Then $f$ is continuous with respect to the Zariski topology on $X$, and hence defines a patch bundle.
Moreover, $f$ is spectral if and only if $U$ is clopen.
\end{example}

The previous example shows that a patch bundle with finite image need not be spectral. By contrast, the following example exhibits a spectral patch bundle with infinite image and illustrates the flexibility in the choice of the indexing Stone space.

\begin{example} \label{chain 1}
Let $R$ be the field ${\mathbb{Q}}$ of rational numbers and consider a strictly increasing
chain of subrings
\[
S_1\subsetneq S_2\subsetneq\cdots
\subsetneq\bigcup_{n\geq 1}S_n=\mathbb Q.
\]
For example, we may take $S_n=\mathbb Z[1/p_1,\ldots,1/p_n]$, where $\{p_n\}$ is an enumeration of the primes. Let $Y= {\mathbb{N}} \cup \{\infty\}$ be the one-point compactification of ${\mathbb{N}}$. Then $Y$ is a Stone space. Define $f:Y\to\Sigma(\mathbb Q)$ by $f(n)=S_n$ for $n\in\mathbb N$ and $f(\infty)=\mathbb Q$. For every $r\in\mathbb Q$, the set $f^{-1}(\mathcal U(r))$ is a clopen tail of $Y$. Hence $f$ is a spectral patch bundle. Moreover, $f$ is a homeomorphism onto its image
$X=\{S_n:n\in\mathbb N\}\cup\{\mathbb Q\}$, endowed with the subspace topology induced by
$\Sigma(\mathbb Q)^{\rm patch}$. Every infinite Stone space admits a continuous
surjection onto $Y$ (this follows, for example, from Stone duality and the proof of
\cite[Lemma~4.2]{BH}). Consequently, the same space of subrings $X$ can be continuously indexed by any infinite Stone space.

The construction also extends to continuous increasing chains of subrings indexed by compact ordinal intervals, although not every infinite Stone space necessarily admits a continuous surjection onto such an interval. Here, continuity means that the ring at each nonzero limit ordinal is the union of the preceding rings.
\end{example}

The next example motivates both why we want the flexibility of working with patch bundles rather than just patch spaces of subrings, and why morphisms of patch bundles are defined in a more nuanced way than simply requiring a commutative diagram of maps.  

\begin{example}\label{ex: zariski inclusion} 
Let $K$ be a field, and let ${\rm Zar}(K)$ denote the Zariski-Riemann space of $K$, that is, the set of valuation rings of $K$ with the subspace Zariski topology inherited from $\Sigma(K)$.
Since ${\rm Zar}(K)$ is a patch space of $\Sigma(K)$, the inclusion map
\[
\iota\colon\; {\rm Zar}(K)^{\rm patch}
\longrightarrow \Sigma(K)
\]
is the tautological patch bundle associated with
${\rm Zar}(K)$. 

Now let $F$ be a field extension of $K$. The contraction map
\[
\phi\colon\; {\rm Zar}(F)^{\rm patch} \longrightarrow  {\rm Zar}(K)^{\rm patch}\colon\; \,V \longmapsto V \cap K
\]
is continuous and surjective. Hence its composition with $\iota$ gives a patch bundle
\[
\iota\circ\phi\colon\; {\rm Zar}(F)^{\rm patch}
\longrightarrow \Sigma(K).
\]
The fiber of $\phi$ over a valuation ring $W$ of $K$
consists precisely of the valuation rings of $F$
extending $W$. Thus the new indexing space
distinguishes different extensions of the same
valuation ring.

Viewing the valuation rings of $K$ as subrings of $F$, we can also express this situation using
a morphism of patch bundles over $F$. Consider the tautological patch bundles
\[
\alpha\colon\; {\rm Zar}(K)^{\rm patch}
\longrightarrow \Sigma(F)
\: \text{ and }\:
\beta\colon\; {\rm Zar}(F)^{\rm patch}
\longrightarrow \Sigma(F).
\]
For each $V\in {\rm Zar}(F)$, we have
\[
\alpha(\phi(V))=V\cap K\subseteq V=\beta(V).
\]
Thus $\phi$ defines a morphism of patch bundles from $\beta$ to $\alpha$. If $F\ne K$, the
inclusion is strict, since $\operatorname{Frac}(V\cap K)=K$ and $\operatorname{Frac}(V)=F$.
In particular, requiring equality instead of inclusion in the definition of a morphism of patch bundles would be too restrictive.
\end{example}

A related version of the following example appears in \cite[Example~3.3(6)]{CosOlb}, in the setting of projective integral schemes. Here we consider more generally integral separated quasi-compact schemes, showing how they can be realized as spaces of local subrings of their function fields and how dominant morphisms naturally induce morphisms of the associated patch bundles.

\begin{example} \label{scheme 1} Fix a field $F$, and let $A$ be a subring of $F$ with quotient field $F$. Let
\[
X(A)=\{A_P:P\in\Spec(A)\}.
\]
Then $X(A)$, viewed as a subspace of $\Sigma(F)$, is homeomorphic to $\Spec(A)$
\cite[Lemma~2.4]{DobbsFedderFontana}.
Moreover, the affine scheme $\Spec(A)$ with its structure sheaf is isomorphic as a locally ringed
space to $X(A)$ equipped with the sheaf whose ring of sections over each nonempty open subset
$U$ is the intersection of the rings in $U$ (and whose ring of sections over the empty set is the zero ring). 
Since $X(A)$ is a spectral space \cite[Theorem~3.2]{SpiritoLoc}, $X(A)^{\rm patch} \rightarrow X(A)$ is the tautological patch bundle associated with $X(A)$ over $F$, where here $X(A)^{\mathrm{patch}}$ denotes the intrinsic patch topology of $X(A)$ as a spectral space.   

This idea also extends to integral separated quasi-compact schemes with function field $F$. Such a scheme $X$ admits a finite affine open cover $X=U_1\cup\cdots\cup U_n$, where $U_i\simeq\Spec(A_i)$ and each $A_i$ is a subring of $F$ with quotient field $F$. Identifying each point $x\in X$ with its local ring $\mathcal O_{X,x}$, we can thus realize $X$ as
\[
X(A_1)\cup\cdots\cup X(A_n)\subseteq\Sigma(F).
\]
Since $X$ is quasi-compact and separated, hence spectral \cite[\href{https://stacks.math.columbia.edu/tag/094L}{Tag 094L}]{stacks-project}, the inclusion $X^{\mathrm{patch}}\rightarrow\Sigma(F)$ is the tautological patch bundle associated with $X$, where $X^{\mathrm{patch}}$ denotes the intrinsic patch topology of $X$. 

Morphisms of patch bundles arise naturally in this setting.  Let $h\colon\; X'\to X$ be a dominant morphism of integral separated quasi-compact schemes with function fields $F'$ and $F$, Since $h$ is dominant, it induces an embedding $F\hookrightarrow F'$. Identifying $F$ with its image in $F'$, we can regard the local rings of both schemes as subrings of $F'$, and consider the tautological patch bundles $$
\alpha\colon\; X^{\mathrm{patch}}\longrightarrow\Sigma(F')\colon\; x\longmapsto\mathcal O_{X,x}\quad\text{and}\quad\beta\colon\; (X')^{\mathrm{patch}}\longrightarrow\Sigma(F')\colon\; x'\longmapsto\mathcal O_{X',x'}. $$
Since $X'$ is quasi-compact and $X$ is quasi-separated, $h$ is quasi-compact and hence spectral. It therefore induces a continuous map $h^{\rm patch}\colon\; (X')^{\mathrm{patch}}\to X^{\mathrm{patch}}$ (see \cite[\href{https://stacks.math.columbia.edu/tag/03GI}{Tag 03GI}]{stacks-project}
and \cite[\href{https://stacks.math.columbia.edu/tag/0G1J}{Tag 0G1J}]{stacks-project}). Moreover, for each $x'\in X'$ the induced homomorphism on local rings is identified with the inclusion $ \alpha(h(x'))=\mathcal O_{X,h(x')} \subseteq \mathcal O_{X',x'}=\beta(x').$ Thus $h^{\mathrm{patch}}$ defines a morphism of patch bundles $\beta\to\alpha$.
\end{example}

\subsection[Patch presheaves]{Patch presheaves}\label{subsect: patch presheaves and bundles}
A {\it patch presheaf} $({{\mathcal R}}, B)$ for a ring $R$ is a collection ${{\mathcal R}}=\{R_e:e \in B\}$ of subrings of $R$ indexed by a Boolean algebra $B$ such that $R=R_0$ (with $0$ the bottom element of $B$), $R_e \cap R_f = R_{e \vee f}$ for all $e,f \in B$, and so also $R_f \subseteq R_e$ whenever $e \leq f$.  

Given a patch presheaf $({{\mathcal R}}, B)$ for $R$, we define its {\it stalks} as the rings \[{{\mathcal R}}_{\m} = \bigcup_{e \not \in {\ff m}} R_e,\] where ${\mathfrak{m}}$ is a maximal ideal in $B$, and denote their collection by \[{X}({{\mathcal R}})=   \{{{\mathcal R}}_{\m}:\m \in \Max(B)\}.\]

Following \cite[Definition 3.10]{CosOlb}, a {\it morphism of patch presheaves for $R$} 
\[({\mathcal R_1}=\{R_e : e\in B_1\}, B_1)\to ({\mathcal R_2}=\{S_f: f\in B_2\}, B_2)\]
consists of a morphism of Boolean algebras $h\colon\; B_2\to B_1$, together with ring inclusions $S_f\subseteq R_{h(f)}$, where $f$ varies in $B_2$. We denote by $\mathsf{PPre}_R$ the category of patch presheaves for $R$, with the morphisms just described.
 
\medskip
Every patch bundle $f\colon\; Y \rightarrow \Sigma(R)$ over $R$ gives rise to a patch presheaf in a natural way: 
\[\mathcal P(f)=(\{R_U:U\in\Clop(Y)\},\Clop(Y)), \: {\mbox{ where }}
R_U=\bigcap_{y\in U}f(y) {\mbox{ and }} R_\emptyset =R,\]
 is a patch presheaf for $R$. 

This construction of a patch presheaf from a patch bundle can be extended to a functor
\begin{equation}\label{eq: functor P}
\mathcal P\colon\; \mathsf{PBun}_R\to \mathsf{PPre}_R 
\end{equation}
as follows.
On objects, $\mathcal P$ sends a patch bundle $f\colon\; Y\to\Sigma(R)$ to $\mathcal{P}(f)$. On morphisms, $\mathcal P$ sends a patch bundle morphism $f_1\to f_2$ represented by $g\colon\; Y_1\to Y_2$ to the morphism of patch presheaves defined by
\[
g^*\colon\; \Clop(Y_2)\longrightarrow\Clop(Y_1)\colon\; V\longmapsto g^{-1}(V),
\]
and the inclusions $R_{2,V}\subseteq R_{1,g^{-1}(V)}$, where $R_{i,U}=\bigcap_{y\in U}f_i(y)$. These inclusions follow from $f_2(g(y))\subseteq f_1(y)$.

Conversely, there is also a natural way   to associate a patch bundle to a patch presheaf for $R$. We define a functor that will turn out to be the  quasi-inverse of ${\mathcal P}$:
\begin{equation}\label{eq: functor P-1}
\mathcal P^{-1}\colon\; \mathsf{PPre}_R\to \mathsf{PBun}_R.
\end{equation}
On objects, $\mathcal P^{-1}$ sends a patch presheaf $(\mathcal R,B)$ to its stalk bundle
\[
\mathcal P^{-1}(\mathcal R,B)=\phi_{\mathcal R}\colon\; \Max(B)\longrightarrow\Sigma(R)\colon\; \m\longmapsto\mathcal R_{\m}.
\]
On morphisms, let $\psi\colon\; (\mathcal R_1,B_1)\to(\mathcal R_2,B_2)$ be defined by $h\colon\; B_2\to B_1$ and inclusions $R_{2,e}\subseteq R_{1,h(e)}$, where $\mathcal R_i=\{R_{i,e}:e\in B_i\}$. Then $\mathcal P^{-1}(\psi)$ is the bundle morphism represented by the Stone-dual map
\[
\Max(B_1)\longrightarrow\Max(B_2)\colon\; \m\longmapsto h^{-1}(\m).
\]
The ring inclusions imply
\[
\phi_{\mathcal R_2}(h^{-1}(\m))
=\bigcup_{e\notin h^{-1}(\m)}R_{2,e}
\subseteq\bigcup_{d\notin\m}R_{1,d}
=\phi_{\mathcal R_1}(\m),
\]
so this map is a morphism of patch bundles.

The assignments in \eqref{eq: functor P} and \eqref{eq: functor P-1} preserve identities and composition by the corresponding properties of inverse image, and hence define functors.

It will be useful in the following to recall from \cite[Lemma 4.2]{CosOlb} that, given a patch bundle $f\colon\; Y \rightarrow \Sigma(R)$, for every $y\in Y$, $f(y)$ is the stalk of the patch presheaf ${\mathcal P}(f)$ on the maximal ideal $\m=\{U\in {\rm Clop}(Y):y \not \in U\}$ of ${\rm Clop}(Y)$ induced by $y$. Namely,

\begin{equation}\label{eq: stalk-reconstruction}
f(y)=\bigcup_{U\ni y} R_U =: {\mathcal P}(f)_y.
\end{equation}

Thus the subrings indexed by $Y$ are the stalks of the patch presheaf of the patch bundle.

\begin{theorem}[{\cite[Theorem 4.5]{CosOlb}}]\label{thm: bundle presheaf equivalence}
The functors $\mathcal P$ and $\mathcal P^{-1}$ are quasi-inverse equivalences between the category of patch bundles over $R$ and the category of patch presheaves for $R$.
\end{theorem}

Here and below, $(-)^{-1}$ denotes a chosen quasi-inverse, rather than a strict inverse.

We consider again the patch bundle examples from Subsection~\ref{subsec: patch bundles} and describe their patch presheaves.

\begin{example} \label{finite 2}
Let $X$ be a finite set of subrings of the ring $R$ and consider the tautological patch bundle $f\colon\; X^{\rm patch} \rightarrow X \subseteq \Sigma(R)$. Then the Boolean algebra $B$ of clopens of $X^{\rm patch}$ is the power set of $X$, and so the patch presheaf associated to $f$ is indexed by
the power set of $X$, and assigns to each subset $E\subseteq X$ the ring $R_E=\bigcap_{T\in E}T$, where $R_\varnothing = R$.
The stalks are the rings belonging to $X$ since each $T\in X$ corresponds to the maximal ideal $\mathfrak m_T=\{E\subseteq X:T\notin E\}$
of $B$.

Now consider, as in Example~\ref{finite 1}, a patch bundle $f\colon Y\to X$, where $X=\{S,T\}$ and $S\subsetneq T$. Let $U=f^{-1}(\{T\})$, which is open in $Y$. The patch presheaf ${\mathcal P}(f)$ has Boolean algebra $B=\operatorname{Clop}(Y)$ and, for each clopen $V$ in $Y$, its
coefficient ring is
\[
R_V=
\begin{cases}
R & \text{if } V=\varnothing,\\
T & \text{if } \varnothing\ne V\subseteq U,\\
S & \text{otherwise}.
\end{cases}
\]
Indeed, if $V$ meets $Y\setminus U$, then the intersection of the rings $f(y)$, for $y\in V$, is $S$. 
%and for each clopen $V$ in $Y$, the coefficient ring $R_V$ is $R$ if $V$ is empty, $T$ if $\varnothing\ne V\subseteq U$ and $S$ otherwise.
Consequently, for a maximal ideal $ \mathfrak m_y=\{V\in{\rm Clop}(Y):y\notin V\}$ of Clop$(Y)$,  the stalk $
\mathcal R_{\mathfrak m_y}$ is $
T$ if $y\in U$ and $S$ otherwise. 
\end{example}

\begin{example} \label{chain 2}
Consider the patch bundle of Example~\ref{chain 1},
$f\colon\; Y=\mathbb N\cup\{\infty\}\rightarrow\Sigma(\mathbb Q),$ where $ f(n)=S_n $ and $f(\infty)=\mathbb Q.$
Its associated patch presheaf is
$$
\mathcal P(f)=(\mathcal R,B), {\mbox{ 
where }}
\mathcal R=\{R_U:U\in B\} {\mbox{ and }}  
B={\rm Clop}(Y).$$ The clopen subsets of $Y$ are precisely the finite subsets of $\mathbb N$ and their complements in $Y$. For $U\in {\rm Clop}(Y)$, since the rings $S_n$ form an increasing chain, we have $R_U = {\mathbb{Q}}$ when $U \cap {\mathbb N}$ is empty, i.e.,  $U= \emptyset$, and $R_U = S_n$, with $n = \min(U \cap {\mathbb N})$, otherwise. 

In particular, $R_{\{n\}}=S_n$, so the stalk at $n$ is $S_n$. On the other hand, the clopen neighborhoods of $\infty$ contain tails of
$\mathbb N$, and hence the stalk at $\infty$ is
\[
\mathcal R_{\mathfrak m_\infty}
=\bigcup_{n\geq 1}S_n=\mathbb Q.
\]
Thus the stalk bundle recovers $f$. Notice that, although the stalk at $\infty$ is $\mathbb Q$, no nonempty clopen subset
has coefficient ring $\mathbb Q$.
\end{example}

\begin{example}\label{ex: zariski inclusion 2} Let $K$ be a field. For the tautological patch bundle $\iota\colon\; {\rm Zar}(K)^{\mathrm{patch}}\to \Sigma(K)$ of Example~\ref{ex: zariski inclusion}, the patch presheaf $\mathcal P(\iota)$ is indexed by
$B=\operatorname{Clop}(\operatorname{Zar}(K)^{\mathrm{patch}})$, with coefficient rings
$R_U=\bigcap_{W\in U}W,$ where  the empty intersection is understood to be $K$. By the stalk identity \eqref{eq: stalk-reconstruction}, the stalk of the presheaf at each $W\in\operatorname{Zar}(K)$ is precisely $W$.

Now let $F$ be a field extension of $K$, and consider the patch bundle $\iota \circ \phi\colon\; {\rm Zar}(F)^{\rm patch} \rightarrow \Sigma(K)$ from Example~\ref{ex: zariski inclusion}, where $\phi(V)=V\cap K$.
Its associated patch presheaf is indexed by $\operatorname{Clop}(\operatorname{Zar}(F)^{\mathrm{patch}})$, with coefficient rings
$R_E=\bigcap_{V\in E}(V\cap K).$ Again by \eqref{eq: stalk-reconstruction}, the stalk at each $V\in\operatorname{Zar}(F)$ is $V\cap K.$
Thus distinct valuation rings of $F$ with the same contraction to $K$ give rise to the same stalk, even though they correspond to distinct
points of the indexing Stone space.
\end{example}

\begin{example}\label{scheme 2}
Let $X$ be an integral separated quasi-compact scheme with function field $F$, and let $f\colon\;X^{\rm patch}\to \Sigma(F)\colon\; x \mapsto \mathcal O_{X,x}$ be its tautological patch bundle as in Example~\ref{scheme 1}. The patch presheaf $\mathcal P(f)$ has indexing Boolean algebra $ B_X=\operatorname{Clop}(X^{\mathrm{patch}})$  and coefficient rings $R_U=\bigcap_{x\in U}\mathcal O_{X,x},$ where $R_\varnothing=F$. By the stalk identity \eqref{eq: stalk-reconstruction}, $\mathcal P(f)_x=\mathcal O_{X,x}$ for each $x\in X$. Every nonempty affine open subset
$U\simeq\Spec(A)$ of $X$ belongs to $B_X$, since it is quasi-compact and open. Moreover,
\[
R_U=\bigcap_{x\in U}\mathcal O_{X,x}
=\Gamma(U,\mathcal O_X)=A.
\]
Now let $h\colon X'\to X$ be a dominant morphism of such schemes, as in Example~\ref{scheme 1}, and identify $F$
with its image in the function field $F'$ of $X'$. Let $\alpha$ and $\beta$ be the tautological patch bundles of $X$ and $X'$,
respectively, viewed over $F'$. The morphism of patch bundles $\beta\to\alpha$ induced by $h^{\rm patch}$ gives, upon applying $\mathcal P$, a
morphism $\mathcal P(\beta)\to\mathcal P(\alpha)$ with underlying Boolean homomorphism $h^{-1}\colon\; B_X\to B_{X'}$. Thus the patch presheaf associated with $X$ recovers its local rings as stalks and its affine coordinate rings as coefficient rings, and the construction is compatible with dominant morphisms of schemes after passing to a common field of functions.
\end{example}

The equivalence $\mathcal P$ of Theorem~\ref{thm: bundle presheaf equivalence} is the first of three equivalences from the category of patch bundles over $R$. The equivalence $\mathcal A$ to the opposite category $\mathsf{PAlg}_R^{\mathrm{op}}$ of patch algebras for $R$ will be established in Section~\ref{sec: patch algebras and patch bundles}, and the  equivalence $\mathcal S$ to the category $\mathsf{PRS}_R$ of patch ringed spaces over $R$ in Section~\ref{sec: patch ringed spaces}:
\[
\begin{tikzcd}[column sep=large, row sep=large]
& \mathsf{PBun}_R
  \arrow[dl, "\mathcal P"']
  \arrow[d, "\mathcal S"]
  \arrow[dr, "\mathcal A"] & \\
\mathsf{PPre}_R & \mathsf{PRS}_R & \mathsf{PAlg}_R^{\mathrm{op}}
\end{tikzcd}
\]

\section[Patch algebras and patch bundles]{Patch algebras and patch bundles}\label{sec: patch algebras and patch bundles}

In \cite{CosOlb}, patch algebras were introduced for patch presheaves and  patch bundles. In this section, we provide a   free-standing notion of a patch algebra   that, {\it a priori}, does not involve these two notions. We then construct in Subsection~\ref{subsec: patch algebras of bundles} a duality between $\mathsf{PBun}_R$ and the category of patch algebras for $R$.
We show that a patch algebra may be viewed concretely as a ring of locally constant $R$-valued functions on a Stone space subject to certain membership conditions. Alternatively, it can be viewed as a subring of a Specker algebra containing a faithful generating Boolean algebra of idempotents. Both descriptions are discussed in this section.

\subsection[Patch algebras and their morphisms]{Patch algebras and their morphisms}

\begin{definition}\label{dfn: patch algebra}
 A {\it Specker presentation (over $R$)} of a ring $A$ consists of a Boolean subalgebra $B_A$ of $\Id(A)$ and an injective ring homomorphism
$$
\iota_A\colon\; A\longrightarrow R[B_A]
$$
whose restriction to $B_A$ is the canonical embedding of $B_A$ into $R[B_A]$. A {\it patch algebra for $R$} is a ring $A$ equipped with a Specker presentation.
\end{definition}

\begin{remark}
It is implicit in Definition~\ref{dfn: patch algebra}  that a patch algebra for the ring $R$ is a triple {$(B_A, A, \iota_A)$, where $B_A$ and $\iota_A$ form a Specker presentation of the ring $A$}. Two Specker presentations of the same ring $A$ therefore determine two distinct patch algebras. Throughout this work, when we refer to a patch algebra $A$,  it is understood that it comes equipped with an associated Boolean algebra $B_A$ and {an embedding $\iota_A\colon\; A\to R[B_A]$}, that is, with a fixed Specker presentation of $A$.

This is equivalent to specifying inclusions $B_A\subseteq A\subseteq S$, where $S$ is a Specker $R$-algebra with faithful generating Boolean algebra $B_A$, since the unique $R$-algebra isomorphism $S\to R[B_A]$ fixing $B_A$ restricts to $\iota_A$ (\cite[Theorem~2.7]{BMMO}), and, conversely, we can identify $A$ with $\iota_A(A)$ and take $S=R[B_A]$. 
\end{remark}

\begin{remark}
Let $A$ be a patch algebra for $R$ with Specker presentation {$(B_A,A,\iota_A)$}. As proved in \cite[Lemma 2.1]{BMMO}, every element  $a \in R[B_A]$ has a unique full orthogonal decomposition with distinct coefficients, and therefore, if $a \neq 0$, it also has a unique orthogonal decomposition with distinct nonzero coefficients. Via $\iota_A$, the elements of $A$ have these decompositions as well; hence the notations $\fod(a)$ and $\od(a)$ introduced in Subsection \ref{subsec: Sepecker algebras} also apply to elements of $A$, and we will make extensive use of them. In calculations, we identify $A$ with $\iota_A(A)\subseteq R[B_A]$; thus an expression $a=\sum_i r_ie_i$ means $\iota_A(a)=\sum_i r_ie_i$ in $R[B_A]$.
\end{remark}

\begin{definition} Let $A$ be a patch algebra for $R$. 
For each $e \in B_A$, define
    $$(A:_R e):=\{r\in R : re\in {\iota_A(A)}\},$$
where the product $re$ is computed in $R[B_A]$. Under the preceding identification, we continue to write $re\in A$ for $re\in\iota_A(A)$. The set $(A:_R e)$ is a subring of $R$
%: it contains $1_R$ because $e\in A$, and it is closed under subtraction and multiplication since $(r-s)e=re-se$, and $(rs)e=(re)(se)$. We call $(A:_R e)$ the {\it coefficient ring of $A$ at $e$}. 
and we call it the {\it coefficient ring of $A$ at $e$}.
\end{definition}

In the following, we will drop the $R$ in the subscript and simply write $(A: e)$ when no confusion can arise.

\begin{lemma}\label{lem: presheaf of a patch algebra}
    Let $A$ be a patch algebra for $R.$ Then for all $e,f \in B_A$, 
   $$ (A: e) \cap (A: f)=(A: e+f-ef).$$
\end{lemma}
\begin{proof}
If $e,f \in B_A$ and $r \in R$ is such that $re, rf \in A$, then $r(e+f-ef) \in A$ as well and thus 
        $$(A: e) \cap (A: f) \subseteq (A: e+f-ef).$$ Vice versa, if $r(e+f-ef) \in A,$ multiplying by $e \in A$ we
        get $r(e+f-ef)e = re \in A,$ so $r \in (A:e).$ An analogous argument shows that $r \in (A:f)$ and
        concludes the proof.
    \end{proof}

For a patch algebra $A$, Lemma~\ref{lem: presheaf of a patch algebra} and the identity $(A:0)=R$ show that 
\[(\{(A:e):e\in B_A\},\:B_A)\]
is a patch presheaf, since $e\vee f=e+f-ef$ in $B_A$.

   \begin{definition}\label{def: morphism of patch algebras}
Let $A_1$ and $A_2$ be patch algebras for $R$. A {\it morphism $\psi\colon\; A_1\to A_2$ of patch algebras} is a ring homomorphism that restricts to a map from $B_{A_1}$ to $B_{A_2}$ and satisfies
\[
\iota_{A_2}\circ\psi=\widetilde{\psi}\circ\iota_{A_1},
\]
where
$
\widetilde{\psi}\colon\; R[B_{A_1}]\longrightarrow R[B_{A_2}]
$
is the unique $R$-algebra homomorphism extending the Boolean algebra homomorphism $\psi|_{B_{A_1}}\colon\; B_{A_1}\to B_{A_2}$ (\cite[Lemma 2.5]{BMMO}). Thus the following diagram commutes:
\[
\begin{tikzcd}[row sep=large,column sep=large]
B_{A_1}
  \arrow[r, hook]
  \arrow[d, "\psi|_{B_{A_1}}"']
& A_1
  \arrow[r, hook, "\iota_{A_1}"]
  \arrow[d, "\psi"]
& R[B_{A_1}]
  \arrow[d, "\widetilde{\psi}"]
\\
B_{A_2}
  \arrow[r, hook]
& A_2
  \arrow[r, hook, "\iota_{A_2}"']
& R[B_{A_2}]
\end{tikzcd}
\]
\end{definition}

After identifying each $A_i$ with its image under $\iota_{A_i}$, the map $\psi$ is the restriction of $\widetilde{\psi}$.
Also, Definition \ref{def: morphism of patch algebras} implies that $(A_1: e)\subseteq (A_2: \psi(e))$, for every $e\in B_{A_1}$.
Together with the induced Boolean homomorphism, these inclusions define a morphism between the associated patch presheaves in the reverse direction.

\begin{remark}
Let us observe that the requirement that $\psi(B_{A_1}) \subseteq B_{A_2}$ is included in the definition of morphism of patch algebras, as it does not follow from the fact that $\psi$ is the restriction of an $R$-algebra homomorphism of Specker algebras. Consider, for example, the case of a decomposable ring $R$, and a finite Boolean algebra $B$ such that $R[B]$ contains a Boolean algebra $B'$ of idempotents different from $B$ such that $R[B]=R[B']$ (such an
example can be found in \cite[Example~3.4]{BMMO}). Consider the identity map $\psi \colon\; R[B] \to R[B'] $; this is an $R$-algebra morphism but
$\psi(B) \nsubseteq B'$. Of course, if we assume that $R$ is indecomposable, the fact that $\psi( B_{A_1}) \subseteq B_{A_2}$ follows immediately from the fact that $\psi$ is a ring homomorphism.
\end{remark}

\begin{remark}
    The restricted map $\psi_{|B_{A_1}} \colon\; B_{A_1} \to B_{A_2}$ is a Boolean algebra morphism; this follows easily from the fact that $\psi$ is a ring homomorphism.
\end{remark}

With the morphisms defined in Definition~\ref{def: morphism of patch algebras}, patch algebras for $R$ form a category, which we denote by $\mathsf{PAlg}_R$. In the next subsection, we establish a category equivalence between $\mathsf{PBun}_R$ and $\mathsf{PAlg}_R^{\mathrm{op}}$.

\subsection[Patch algebras associated to patch bundles]{Patch algebras associated to patch bundles}\label{subsec: patch algebras of bundles}
In \cite[Section 6]{CosOlb}, the first two authors defined the {\it patch algebra of a patch presheaf $({{\mathcal R}}, B)$ for $R$}, denoted by
${{\mathcal R}}[B]$, as the subring of the Specker $R$-algebra $R[B]$ defined by 
\smallskip
\begin{center}
${{\mathcal R}}[B] = \{t_1{e_1} + \cdots + t_n {e_n} : \; n>0,\ e_i \in B$ and $t_i \in R_{e_i}$ for every $i\in \{1,\ldots, n\}\}$.
\end{center}

In light of the equivalence between the category of patch bundles and the category of patch presheaves for a ring $R$ established in Theorem~\ref{thm: bundle presheaf equivalence}, we can define in a natural way the {\it patch algebra of a patch bundle $f\colon\; Y \rightarrow X$ over $R$} as the patch algebra of the patch presheaf for $R$ given by
   $$({{\mathcal R}}, B)={\mathcal P}(f)=(\{R_U : U\in {\rm Clop}(Y)\}, {\rm Clop}(Y)).$$
 We denote this algebra by $A(f)$, as defined in \cite[Definition 6.7]{CosOlb}. Here, $B_{A(f)} = {\rm Clop}(Y)$ and $\iota_{A(f)}\colon\; A(f)\hookrightarrow R[B_{A(f)}]$ is the inclusion; so $A(f)={\mathcal R}[B_{A(f)}]$ is a patch algebra with Specker presentation $$B_{A(f)} \subseteq A(f)  \subseteq R[B_{A(f)}].$$

An important feature of this construction is that the stalks of $A(f)$ recover the rings in the image of the patch bundle:

 \begin{proposition}[{\cite[Lemma~6.8 and its proof]{CosOlb}}]
\label{prop: stalks as factor rings}
Let $f\colon\; Y\to\Sigma(R)$ be a patch bundle. For each $y\in Y$, let $\mathfrak m_y=\{U\in\Clop(Y):y\notin U\}$ be its corresponding maximal ideal of $B_{A(f)}$.
There is a canonical surjective ring homomorphism
$A(f)\to f(y)$ with kernel $\mathfrak m_yA(f)$, the ideal
generated by the idempotents in $\mathfrak m_y$.
Consequently,
\[
A(f)/\mathfrak m_yA(f)\cong\mathcal P(f)_y=f(y).
\]
\end{proposition}

We next give a different interpretation of $A(f)$ that will be useful in the next sections. We do so by representing $A(f)$ as a ring of
continuous functions.

\medskip

For a Stone space $Y$, the {\it canonical Specker representation} of $Y$ is the $R$-algebra isomorphism  (see \cite[proof of Theorem~2.7]{BMMO})
$$
R[\Clop(Y)]\longrightarrow C(Y,R_{\mathrm{disc}})\colon\; e_U\longmapsto 1_U,
$$
where $e_U$ is the idempotent corresponding to the clopen subset $U$. If a subring $S\subseteq C(Y,R_{\mathrm{disc}})$
contains all characteristic functions of clopens, its {\it canonical Specker presentation} has Boolean algebra
$B_S=\{1_U:U\in\Clop(Y)\}$ and embedding
$$
S\mathrel{\scalebox{1.6}[1]{$\hookrightarrow$}} C(Y,R_{\mathrm{disc}})
\overset{\sim}{\longrightarrow}R[B_S],
$$
where the second map is the canonical isomorphism fixing $B_S$.

\begin{theorem}\label{C rep}
Let $f\colon\; Y\to\Sigma(R)$ be a patch bundle over $R$. The canonical
Specker representation of $Y$ restricts to an isomorphism of patch algebras
$$
A(f)\cong
\{\sigma\in C(Y,R_{\mathrm{disc}}):
  \sigma(y)\in f(y)\text{ for every }y\in Y\},
$$
where the ring on the right is endowed with its canonical Specker presentation. Moreover, if a patch bundle morphism
$f_1\to f_2$ is represented by $g\colon Y_1\to Y_2$, then the induced patch algebra morphism $A(f_2)\to A(f_1)$ corresponds, under these identifications,
to precomposition with $g$.
\end{theorem}

\begin{proof}
Let $B=\Clop(Y)$ and, for each $e\in B$, write $R_e=\bigcap_{y\in U_e}f(y)$.
Consider the canonical Specker representation
\[
\phi\colon R[B]\longrightarrow C(Y,R_{\mathrm{disc}}),
\]
which sends each $e\in B$ to the characteristic
function $1_{U_e}$.

Let $a\in R[B]$, with full orthogonal decomposition $\fod(a)=\sum_i a_i e_i$. Since the sets $U_{e_i}$ form a pairwise disjoint
clopen cover of $Y$, the function $\phi(a)$ takes the constant value $a_i$ on $U_{e_i}$. Hence
$$
\phi(a)(y)\in f(y)\text{ for all }y\in Y
\:\:\Longleftrightarrow\:\:
a_i\in\bigcap_{y\in U_{e_i}}f(y)=R_{e_i}\text{ for all }i
\:\: \Longleftrightarrow\:\:a\in A(f).
$$
Thus $\phi$ restricts to the required ring isomorphism. Moreover, it identifies the specified Boolean algebra $B$ with the Boolean
algebra of characteristic functions of clopen subsets of $Y$, and is induced by an isomorphism of the ambient Specker algebras. Hence its
restriction is an isomorphism of patch algebras.

Finally, for a bundle morphism $f_1\to f_2$ represented by $g\colon\; Y_1\to Y_2$, the inclusion $f_2(g(y))\subseteq f_1(y)$ shows that composition with $g$ sends the function ring for $f_2$ into that for $f_1$. Let $h=g^*\colon\; \Clop(Y_2)\to\Clop(Y_1)$. Then
$1_{U_{h(e)}}=1_{U_e}\circ g.$ Thus composition with $g$ corresponds to the restriction of the Specker homomorphism induced by $h$, and is a morphism of patch algebras.
\end{proof}

We now define the functors underlying  the equivalence between patch bundles and patch algebras.

\begin{definition}\label{rem: PAlg PBundle}
Define functors
\[
\mathcal A\colon\; \mathsf{PBun}_R\longrightarrow\mathsf{PAlg}_R^{\mathrm{op}}
\: \text{ and } \:
\mathcal A^{-1}\colon\; \mathsf{PAlg}_R^{\mathrm{op}}\longrightarrow\mathsf{PBun}_R
\]
as follows.  On objects, $\mathcal A$ sends a patch bundle $f$ to $A(f)$. If a morphism $f_1\to f_2$ of patch bundles is represented by $g\colon\; Y_1\to Y_2$, then $\mathcal A(g)$ corresponds to the patch algebra morphism $A(f_2)\to A(f_1)$ induced by
\[
g^*\colon\; \Clop(Y_2)\longrightarrow\Clop(Y_1).
\]
The functor $\mathcal A^{-1}$ sends a patch algebra $A$ to the stalk bundle
\[
f_A\colon\; \Max(B_A)\longrightarrow\Sigma(R) :
\mathfrak m\longmapsto\bigcup_{e\notin\mathfrak m}(A:e).
\]
If $\psi\colon\; A_1\to A_2$ is a morphism of patch algebras, then $\mathcal A^{-1}(\psi^{\mathrm{op}})$ is represented by the Stone-dual map
\[
(\psi|_{B_{A_1}})^*\colon\; \Max(B_{A_2})\longrightarrow\Max(B_{A_1}).
\]
\end{definition}

To see that the assignments in Definition~\ref{rem: PAlg PBundle} do indeed define functors, observe that by Theorem~\ref{C rep},  $\mathcal A(g)$ is well defined.
  By Lemma~\ref{lem: presheaf of a patch algebra}, the coefficient rings of $A$ form a patch presheaf, so the functor $\mathcal P^{-1}$ in \eqref{eq: functor P-1} shows that $f_A$ is a patch bundle.
If $\psi\colon\; A_1\to A_2$ is a morphism of patch algebras, then  Definition~\ref{def: morphism of patch algebras} gives
$$
(A_1:e)\subseteq(A_2:\psi(e)), {\mbox{ for all }} e\in B_{A_1}.
$$
 Writing $h=\psi|_{B_{A_1}}$, these inclusions imply
\[
f_{A_1}(h^{-1}(\mathfrak m))\subseteq f_{A_2}(\mathfrak m)
\:\text{ for every  }\mathfrak m\in\Max(B_{A_2}),
\]
so the Stone-dual map defines the required bundle morphism.
Identities and compositions are preserved by construction.

\begin{theorem}\label{cor: PAlg PBundle}
The functors $\mathcal A$ and $\mathcal A^{-1}$ are quasi-inverse equivalences. In particular, the categories $\mathsf{PBun}_R$ and $\mathsf{PAlg}_R^{\mathrm{op}}$ are equivalent.
\end{theorem}
\begin{proof}
We can identify each patch algebra $A$ with its image under $\iota_A$ in $R[B_A]$. We claim  
\begin{equation}\label{eq: algebra sums}
A=\left\{\sum_{i=1}^{n}r_i e_i:
\; e_i\in B_A, \ r_i\in(A:e_i)\right\}.
\tag{*}
\end{equation}
Every sum on the right belongs to $A$ by the definition of $(A:e_i)$. Conversely, if $a\in A$, write fod$(a)=\sum_i r_i e_i$. Since $r_i e_i=ae_i\in A$, each $r_i$ belongs to $(A:e_i)$. Thus the coefficient rings of $A$ determine $A$ inside its Specker algebra.

Now let $f\colon\; Y\to\Sigma(R)$ be a patch bundle, with associated coefficient rings $R_U$. By Theorem~\ref{C rep}, we have
$(A(f):U)=R_U$ for every $U\in\Clop(Y)$. Consequently, under the canonical Stone homeomorphism $Y\cong\Max(\Clop(Y))$, the bundle $\mathcal A^{-1}\mathcal A(f)$ assigns to each $y\in Y$ the ring $$
\bigcup_{U\ni y}(A(f):U) =\bigcup_{U\ni y}R_U
=f(y),$$
where the last equality is the stalk identity \eqref{eq: stalk-reconstruction}. Thus the canonical Stone homeomorphism induces an isomorphism of patch bundles $f\cong\mathcal A^{-1}\mathcal A(f)$.

Conversely, let $A$ be a patch algebra, and set $f_A=\mathcal A^{-1}(A)$. By Theorem~\ref{thm: bundle presheaf equivalence}, the patch presheaf $\mathcal P(f_A)$ is canonically isomorphic to the coefficient presheaf of $A$. In particular, for each $e\in B_A$ we have $$
\bigcap_{\mathfrak m\in U_e}f_A(\mathfrak m)=(A:e).$$ Therefore, $A(f_A)$ consists of the finite sums in \eqref{eq: algebra sums}, so $\iota_A$ induces an isomorphism of patch algebras $A\cong A(f_A)$.

Finally, the isomorphisms $f\cong\mathcal A^{-1}\mathcal A(f)$ and $A\cong A(f_A)$ are natural in $f$ and $A$, respectively. On the bundle side, this follows from the naturality of the Stone homeomorphisms. On the algebra side, it follows from the defining identity $\widetilde\psi\circ\iota_{A_1}
=\iota_{A_2}\circ\psi$ for a morphism of patch algebras $\psi\colon\; A_1\to A_2$. Thus both composites are naturally isomorphic to the corresponding identity functors.
\end{proof}

We conclude this section by revisiting the running examples from Section~\ref{sec: bundles and presheaves} and describing their associated patch algebras.

\begin{example}\label{finite 3}
Let $X$ be a finite collection of subrings of $R$, and consider the tautological patch bundle $f\colon\; X^{\mathrm{patch}}\to X$ of
Example~\ref{finite 1}. As in Example~\ref{finite 2}, its associated patch presheaf has Boolean algebra $B=\mathcal P(X)$. For each $T\in X$, let $e_T$ denote the idempotent corresponding to the singleton $\{T\}$. These idempotents are pairwise orthogonal and sum to $1$, and their coefficient rings are $R_{\{T\}}=T$. Consequently, every element of $A(f)$ has a
unique expression $\sum_{T\in X}a_Te_T,$ $a_T\in T.$ It follows that $A(f)$ is isomorphic to $\prod_{T\in X}T$.
\end{example}

\begin{example}\label{chain 3}
Consider the chain of rings $(S_n)$ of Example~\ref{chain 1} and the patch bundle $f\colon Y=\mathbb N\cup\{\infty\}\to\Sigma(\mathbb Q)$, where
$f(n)=S_n$ for $n\in\mathbb N$ and $f(\infty)=\mathbb Q$.
Since $\mathbb Q$ has the discrete topology, a continuous function $\sigma\colon \; \mathbb N \cup \{\infty\}\to \mathbb Q$ is continuous at $\infty$ if and only if $\sigma(n)=\sigma(\infty)$ for all sufficiently large $n$. Thus restriction to $\mathbb N$ determines $\sigma$ uniquely, and every eventually constant sequence $(s_n)$ with $s_n\in S_n$ extends continuously to $Y$ by assigning its eventual constant value to $\infty$.
Therefore, by Theorem~\ref{C rep}, restriction to $\mathbb N$ identifies the patch algebra $A=A(f)$ with the ring
$$\left\{ (s_n)\in\prod_{n\ge1}S_n: (s_n)\text{ is eventually constant} \right\},$$
with coordinatewise operations.  

Evaluation at $\infty$ defines a homomorphism $\operatorname{ev}_\infty\colon\;  A\to\mathbb Q$ that sends each eventually constant sequence to its eventual constant value, which is surjective since $\bigcup_{n\geq1}S_n=\mathbb Q$.
Its kernel is the maximal ideal $P_\infty = \bigoplus_{n\ge1}S_n$ where the direct sum is identified with the ideal of finitely supported sequences.

Every subring of $\mathbb Q$ is a localization of $\mathbb Z$, hence a Pr\"ufer domain, so $A$ is semihereditary by \cite[Corollary~7.2]{CosOlb}. Moreover, by applying  \cite[Theorem~6.9 and Proposition~7.4]{CosOlb} we can identify the minimal spectrum $\Min (A)$ of $A$ with the one-point compactification $\mathbb N\cup\{\infty\}$ of ${\mathbb N}$. Its minimal prime ideals are precisely
$P_\infty$ and $P_n=\{(s_i)\in A:s_n=0\}$, for all $\geq1$, with corresponding quotients
\[
A/P_n\cong S_n \: {\mbox{ and }} \: A/P_\infty\cong\mathbb Q.
\]

In particular, $P_\infty$ is both a minimal prime ideal and a maximal ideal, and is the unique nonisolated point of
$\Min(A)$. The ring $A$ is not Noetherian, since its minimal spectrum is infinite.
\end{example}

\begin{example}
In the notation of Example~\ref{ex: zariski inclusion}, let $A=A(\iota)$ be the patch algebra of the tautological patch bundle $\iota\colon\; {\rm Zar}(K)^{\mathrm{patch}}\to \Sigma(K)$. By \cite[Example~7.7]{CosOlb}, the Pierce, minimal, and maximal spectra of $A$ are all homeomorphic to $\operatorname{Zar}(K)^{\mathrm{patch}}$. Under the representation of Theorem~\ref{C rep},
the minimal prime and maximal ideals corresponding
to $W\in\operatorname{Zar}(K)$ are
\[
P_W=\{\sigma\in A:\sigma(W)=0\} \: {\mbox{ and }} \:
M_W=\{\sigma\in A:\sigma(W)\in\mathfrak m_W\},
\]
where $\mathfrak m_W$ is the maximal ideal of $W$.
By Proposition~\ref{prop: stalks as factor rings},
\[
A/P_W\cong A_{M_W}\cong W.
\]
Thus $A$ realizes the valuation rings of $K$ both as quotients by its minimal prime ideals and as localizations at its maximal ideals. 
Moreover, $A$ is clean (i.e., every element is the sum of an idempotent and a unit) by \cite[Proposition~7.5]{CosOlb}, and it is arithmetical because its localizations at maximal ideals are valuation domains. 

Now let $F$ be a field extension of $K$, and consider the patch algebra $A'=A(\iota\circ\phi)$, where $\phi(V)=V\cap K$. Applying \cite[Theorem~6.9 and Propositions~7.4--7.5]{CosOlb},
the minimal prime ideals and maximal ideals of $A'$ are indexed by $V\in\operatorname{Zar}(F)$, with $A'/P'_V\cong A'_{M'_V}\cong V\cap K$ (Proposition~\ref{prop: stalks as factor rings}).
In particular, distinct valuation rings of $F$ with the same contraction to $K$ correspond to distinct minimal prime ideals and distinct maximal ideals of $A'$, even though their associated quotients and localizations are isomorphic to the same valuation ring of $K$. These points remain distinguished by the
idempotents corresponding to clopen subsets of $\operatorname{Zar}(F)^{\mathrm{patch}}$.
\end{example}

\begin{example} 
Let $X$ be an integral separated quasi-compact scheme with function field $F$, and let $A_X=A(f)$, where
$f\colon X^{\mathrm{patch}}\to\Sigma(F)$ is the tautological patch bundle of Example~\ref{scheme 1}, sending each $x\in X$ to $\mathcal O_{X,x}$.
By Theorem~\ref{C rep}, we may identify $A_X$ with the ring
\[
\left\{
\sigma\in C(X^{\mathrm{patch}},F_{\mathrm{disc}}):
\sigma(x)\in\mathcal O_{X,x}
\text{ for every }x\in X
\right\}.
\]
By Theorem~6.9, Propositions~7.4--7.5, and Corollary~7.6 in \cite{CosOlb}, $A_X$ is clean, and its Pierce, minimal, and maximal spectra are homeomorphic to $X^{\mathrm{patch}}$.
More precisely, for each $x\in X$, evaluation at $x$ determines a minimal
prime ideal $P_x=\{\sigma\in A_X:\sigma(x)=0\}$ and a maximal ideal $M_x=\{\sigma\in A_X:\sigma(x)\in\mathfrak m_{X,x}\}$. By Proposition~\ref{prop: stalks as factor rings},
\[
A_X/P_x\cong (A_X)_{M_x}
\cong\mathcal O_{X,x}.
\]
Thus every local ring of $X$ can be recovered from $A_X$, both as a quotient by a minimal prime ideal and as a localization at a maximal ideal.

Now suppose $h\colon\; X'\to X$ is a dominant morphism of such schemes, and assume $F $ is a subfield of the function field $F'$ of $X'$.
Viewing both tautological patch bundles as over $F'$, Theorem~\ref{C rep} yields an induced morphism of patch algebras $h^*\colon\; A_X\rightarrow A_{X'}\colon\; \sigma\mapsto\sigma\circ h.$ Since the inclusion maps $\mathcal O_{X,h(x')}\subseteq\mathcal O_{X',x'}$ are local, $(h^*)^{-1}(M_{x'})=M_{h(x')}$ for every $x'\in X'$. Thus contraction along $h^*$ defines a map $\Max(A_{X'})\to\Max(A_X)$ which, under the canonical homeomorphisms with $(X')^{\mathrm{patch}}$ and $X^{\mathrm{patch}}$, is the map induced by $h$. 

Thus the construction of $A_X$ recovers the local rings and the underlying patch space of $X$, compatibly with dominant morphisms of schemes.
\end{example}

\section{Morphisms of patch algebras and their kernels}\label{sec: morphisms and kernels}

In this section, we investigate morphisms of patch algebras for a ring $R$. In particular, we focus on their kernels, factorizations,
and the geometric characterization of injectivity and surjectivity.

\begin{proposition} \label{prop: surj, inj, bij 1}  Every morphism  of patch algebras for $R$ factors as a surjective morphism followed by an injective morphism. 
\end{proposition}

\begin{proof} 
Let $\psi\colon\; A_1 \rightarrow A_2$ be a morphism of patch algebras, with Specker presentations $$B_i \subseteq A_i \stackrel{\iota_i}{\rightarrow} R[B_i], \:\: i =1,2.$$
 It suffices to show that   $\psi(A_1)$ is a patch algebra with Specker presentation 
 \begin{equation}\label{eq: tag 5.1}
 \psi(B_{1}) \subseteq \psi(A_1) \stackrel{\iota_2}\longrightarrow R[\psi(B_{1})].\    
 \end{equation}
If $ a= \sum_i r_ie_i \in A_1$, where $r_i \in (A_1:e_i)$, then $\psi(a) = \sum_i r_i\psi(e_i)$. Since 
$r_i e_i \in A_1$, we have $r_i\psi(e_i) \in \psi(A_1)$, 
so $r_i \in (\psi(A_1):\psi(e_i))$. Thus \eqref{eq: tag 5.1} is a Specker presentation of $\psi(A_1)$, and so $\psi(A_1)$ is a patch algebra.  
\end{proof}

The kernel of a morphism of patch algebras for $R$ has a particularly transparent form: 

\begin{lemma} \label{lem: surj, inj, bij 2} Let $\psi\colon\; A_1 \rightarrow A_2$ be a morphism of patch algebras, write $B_i=B_{A_i}$ for $i=1,2$, and let ${\mathfrak I} = \ker \psi|_{B_1}$. Then 
$$\ker \psi = {\mathfrak I}A_1=\{\sum r_ie_i:e_i \in {\mathfrak I},r_i \in (A_1:e_i)\}.$$ 
\end{lemma}

\begin{proof} We can assume each $A_i$ has  Specker presentation $B_i \subseteq A_i \:{\hookrightarrow}\: R[B_i]$. Clearly, ${\mathfrak I} \subseteq \ker \psi$. Suppose $a \in \ker \psi$, and write fod$(a) = \sum_ir_ie_i$. Then $0=\psi(a) = \sum_i r_i\psi(e_i)$. Since the $e_i$ are pairwise orthogonal, $0 = r_i\psi(e_i)$ for all $i$. Since $R[\psi(B_1)]$ is a Specker algebra, the $\psi(e_i)$ are faithful and so either $r_i =0$ or $e_i \in {\mathfrak I}$, so we conclude that $a \in {\mathfrak I}A_1$.  Thus $\ker \psi = {\mathfrak I}A_1$.  
For the second equality, suppose $e \in {\mathfrak I}$ and $a \in A_1$. Write fod$(a) = \sum r_ie_i$. Then $ea=\sum_i r_i ee_i$, with $r_i\in(A_1:e_i)\subseteq(A_1:ee_i)$ and $ee_i\in\mathfrak I$.
Thus every element of $\mathfrak I A_1$ is a sum of the desired form.
Conversely, if $r\in(A_1:e)$ and $e\in\mathfrak I$, then
$re=(re)e\in\mathfrak I A_1$.
This proves the second equality.  
\end{proof}

\begin{theorem}\label{thm: kernels R}
    Let $A$ be a patch algebra for a ring $R$ with Specker presentation $B_A \subseteq A \xrightarrow{\iota_A} R[B]$. 
    Then the ideals ${\mathfrak I}A$ of $A$, as $\mathfrak I$ ranges among the ideals of the Boolean algebra
    $B_A\subseteq A$, are precisely the kernels of morphisms of patch algebras for $R$ with source $A$. 
\end{theorem}

\begin{proof}
 Let $\mathfrak I$ be an ideal of $B=B_A$. Then $B' := B/\mathfrak I$ is a Boolean algebra, and there exists a morphism of Boolean algebras $\pi \colon\; B \to B'$
such that $\ker \pi = \mathfrak I.$ Extend $\pi$ to an $R$-algebra homomorphism $ \Phi \colon\; R[B] \to R[B']$. {Set $\psi=\Phi\circ\iota_A$ and $A'=\psi(A)$.}
%{ Let $\phi=\Phi_{|A}$ and $A'=\phi(A)$.} 
We claim that $A'$ is a patch algebra for $R$ and $\psi \colon\; A \to 
A'$ is a morphism of patch algebras such that $\ker \psi = {\mathfrak I}A.$
Since $B' = \pi(B) = \psi(B)$ and $B \subseteq A,$ we have that $B' \subseteq A'\subseteq R[B'].$ Thus, $A'$ is a patch algebra for $R$ with Specker presentation {given by $B'$ and the inclusion $\iota_{A'}\colon\; A'\hookrightarrow R[B']$}. 
Moreover, {$\iota_{A'}\circ\psi=\Phi\circ\iota_A$} by construction, and $B' = \psi(B)$ by definition, so $\psi$ is a morphism of patch algebras for $R$.
Thus, by Lemma~\ref{lem: surj, inj, bij 2}, $\ker \psi = (\ker \psi_{|B})A = (\ker \pi)A = {\mathfrak I}A,$ which proves the forward implication. The converse follows from Lemma~\ref{lem: surj, inj, bij 2}.
\end{proof}

\begin{corollary}\label{cor: kernel correspondence}
Let $A$ be a patch algebra for $R$. The assignment $\mathfrak I\mapsto\mathfrak I A$ is a bijection between the ideals of $B_A$ and the kernels of morphisms of patch algebras for $R$ with source $A$, with inverse $K\mapsto K\cap B_A$.
\end{corollary}

\begin{proof}
By Theorem~\ref{thm: kernels R}, every such kernel has the form $\mathfrak I A$. For each ideal $\mathfrak I$ of $B_A$, the  proof of the theorem gives a morphism $\psi$ with $\ker\psi=\mathfrak I A$ and $\ker(\psi|_{B_A})=\mathfrak I$. Hence $$ \mathfrak I A\cap B_A =\ker\psi\cap B_A =\ker(\psi|_{B_A}) =\mathfrak I, $$
and the corollary follows.
\end{proof}

\medskip

A morphism  $\psi\colon\; A_1 \rightarrow A_2$ of patch algebras is a {ring homomorphism that {preserves the specified Boolean algebras and is compatible with the embeddings $\iota_{A_1}$ and $\iota_{A_2}$} as established in Definition \ref{def: morphism of patch algebras}.} As such, we can consider injectivity and surjectivity of this morphism when considered as a map between the rings $A_1$ and $A_2$.

\begin{theorem}\label{thm: inj surj iso}
Let $A_1$ and $A_2$ be patch algebras for $R$, and let $$f_1\colon\; Y_1 \to \Sigma(R) \: {\mbox{ and }} \: f_2\colon\; Y_2 \to \Sigma(R)$$  be the patch bundles over $R$ corresponding to $A_1$ and $A_2$, respectively. Suppose 
 $\psi\colon\; A_1 \rightarrow A_2$ is a morphism of patch algebras for $R$, and let $g\colon\; Y_2 \rightarrow Y_1$ be the continuous map defining the corresponding morphism of patch bundles.  
   \begin{enumerate}[label=\textup{(\arabic{*})}, mode=unboxed]
   \item\label{thm: inj} The map   $\psi$ is injective if and only if  $g$ is surjective. 
   \item\label{thm: surj} The map  $\psi$ is surjective if and only if $g$ is injective and $f_1 \circ g = f_2$.
   \item\label{thm: iso} The map $\psi$ is an isomorphism if and only if $g$ is a homeomorphism and $f_1 \circ g = f_2$.
   \end{enumerate}
\end{theorem}

\begin{proof} Let $B_i=B_{A_i}$ for $i=1,2$.

\medskip

  \ref{thm: inj} By Lemma \ref{lem: surj, inj, bij 2} and Stone duality, respectively, we have 
  \[\ker \psi =\{0\} \: \Longleftrightarrow \: \ker \psi|_{B_1} =\{0\} \: \Longleftrightarrow \: g \text{ is surjective.}\]
\ref{thm: surj} By Theorem~\ref{C rep}, we can assume  $$A_i=\{\sigma\in C(Y_i,R_{\mathrm{disc}}): \sigma(y)\in f_i(y)\mbox{ for every }y\in Y_i\}$$ and  $\psi(\sigma)=\sigma\circ g$ for all $\sigma \in A_1$.  For each $y\in Y_i$, the map $A_i \rightarrow f_i(y)\colon\; \sigma \mapsto \sigma(y)$  is surjective  since  for each  $r\in f_i(y)$, the nonempty open set
$$\{z\in Y_i:r\in f_i(z)\} 
= f_i^{-1}(\mathcal{U}(r))$$ contains a clopen neighborhood $U$ of $y$, and $r1_U\in A_i$ has value $r$ at $y$.

Suppose first that $\psi$ is surjective. If $g(y)=g(y')$, for all $\sigma \in A_1$, then $$\psi(\sigma)(y) = \sigma(g(y)) = \sigma(g(y')) = \psi(\sigma)(y').$$   Since the characteristic functions of clopen subsets of $Y_2$ belong to $A_2$ and separate points, surjectivity forces $y=y'$: if $y \neq y'$ then we can pick a clopen subset $U$ of $Y_2$ with $y \in U$, $y' \notin U$ and $\sigma \in A_1$ with $\psi(\sigma)= 1_U;$ then
$1=\psi(\sigma)(y) = \psi(\sigma)(y') = 0$ is a contradiction.

Thus $g$ is injective.
Moreover, $f_1\circ g=f_2$ since for every $y\in Y_2$, $$f_2(y) =\{\tau(y):\tau\in A_2\} =\{\sigma(g(y)):\sigma\in A_1\}
=f_1(g(y)).$$

Conversely, suppose that $g$ is injective and $f_1\circ g=f_2$.  Every element of $A_2$ is a finite sum of functions $r1_U$, where $U\subseteq Y_2$ is clopen and $r\in f_2(y)$ for every $y\in U$.
Thus, to verify $\psi$ is surjective, it suffices to show the function $r1_U$ is in the image of $\psi$.
Notice that for such $U$ and $r$ we have $U \subseteq f_2^{-1}(\mathcal{U}(r)) = g^{-1}(f_1^{-1}( \cal U(r)))$, so $g(U) \subseteq f^{-1}_1(\cal U (r))$. Since the latter
set is open and $g(U)$ is compact, we can find a clopen subset $V_1$ of $Y_1$ with $g(U) \subseteq V_1 \subseteq f^{-1}_1(\cal U (r))$.
Injectivity of $g$ and Stone duality imply that the map $g^* \colon\; \Clop(Y_1) \to \Clop(Y_2)$ is surjective, so there is a $V_2 \in \Clop(Y_1)$
with $g^*(V_2)=g^{-1}(V_2)=U$. Taking $V = V_1 \cap V_2 \subseteq f^{-1}_1(\cal U (r))$ we get $g^{-1}(V) = U$ and $r \in f_1(y)$ for all $y \in V$, so that $r1_V \in A_1$ and $\psi(r1_V) =r1_U.$ Thus $\psi$ is surjective.

\medskip

\ref{thm: iso} The morphism of patch algebras $\psi\colon\; A_1\to A_2$ is a ring isomorphism if and only if it is both injective and surjective. In light of items \ref{thm: inj} and \ref{thm: surj}, this happens if and only if $g$ is bijective and $f_1 \circ g = f_2$. Since a continuous bijection between Stone spaces is a homeomorphism, the conclusion follows. 
\end{proof}

\begin{remark}\label{rem: geometric factorization}
In the notation of Theorem~\ref{thm: inj surj iso}, the factorization of $\psi$ through $\psi(A_1)$ in Proposition~\ref{prop: surj, inj, bij 1} corresponds, under the bundle-algebra equivalence, to the factorization of $g$ through $g(Y_2)$, as shown in the following diagrams.
    \[\begin{tikzcd}
	{A_1} && {A_2} \\
	& {\psi(A_1)}
	\arrow["\psi", from=1-1, to=1-3]
	\arrow["\phi"', two heads, from=1-1, to=2-2]
	\arrow["\eta"', hook, from=2-2, to=1-3]
\end{tikzcd}\qquad \qquad\begin{tikzcd}
	 {Y_1} && {Y_2} \\
	& {g(Y_2)}
	\arrow["g"', from=1-3, to=1-1]
	\arrow["{g'}", two heads, from=1-3, to=2-2]
	\arrow["\iota", hook', from=2-2, to=1-1]
\end{tikzcd}.\]
Here $\phi$ and $g'$ are the surjective maps obtained by viewing $\psi$ and $g$ as maps onto their respective images, while $\eta$ and $\iota$ are inclusions. The intermediate patch algebra $\psi(A_1)$ corresponds to the restricted bundle $f_1|_{g(Y_2)}$.
\end{remark}

We now apply the surjectivity criterion in Theorem~\ref{thm: inj surj iso}\ref{thm: surj} to restrictions of patch
bundles to closed subspaces.

\begin{corollary}\label{cor: inclusion of Stone spaces}
    Let $A$ be a patch algebra for $R$, with associated patch bundle
$f\colon\; Y\to\Sigma(R)$. For every closed subspace $Y' \subseteq Y$ there is a 
 surjective morphism of patch algebras \[A\longrightarrow A(f|_{Y'}).\]
\end{corollary}
\begin{proof}
Since $Y'$ is closed in the Stone space $Y$, it is a Stone space.
The inclusion $\;\iota\colon\; Y'\hookrightarrow Y$ defines a morphism of patch
bundles $f|_{Y'}\to f$ and satisfies $f\circ\iota=f|_{Y'}$.
Then the conclusion follows from
Theorem~\ref{thm: inj surj iso}\ref{thm: surj}.
\end{proof}

In particular, let $A$ be a patch algebra for $R$ corresponding to the tautological patch bundle of a patch space $X\subseteq\Sigma(R)$. For every patch closed subspace $X'\subseteq X$, there is a
surjective morphism $A\to A'$, where $A'$ is the patch algebra associated to the tautologial patch bundle of $X'$.

\medskip

For a surjective morphism $\psi\colon\; A\to A'$, Theorem~\ref{thm: inj surj iso}\ref{thm: surj} shows that the dual map
$$(\psi|_{B_A})^*\colon\; \Max(B_{A'})\longrightarrow\Max(B_A)$$ is a homeomorphism onto
its closed image. The next corollary identifies this image in terms
of the Boolean kernel.

\begin{corollary}\label{cor: Boolean kernel closed subspace}
Let $\psi\colon\; A_1\to A_2$ be a surjective morphism of patch algebras for $R$,
write $B_i=B_{A_i}$ for $i=1,2$, and let
$\mathfrak I=\ker \psi|_{B_1}$.
Then $(\psi|_{B_1})^*$ induces a homeomorphism
\[
\Max(B_2)\xrightarrow{\ \sim\ }V(\mathfrak I)
=\{\mathfrak m\in\Max(B_1):\mathfrak I\subseteq\mathfrak m\},
\]
where $V(\mathfrak I)$ is a closed subspace of $\Max(B_1)$.
\end{corollary}

\begin{proof}
By Theorem~\ref{thm: inj surj iso}\ref{thm: surj} and Stone duality,
$\psi|_{B_1}\colon\; B_1\to B_2$ is surjective. Its kernel is $\mathfrak I$, so
$B_1/\mathfrak I\cong B_2$. Stone duality therefore identifies
$(\psi|_{B_1})^*$ with the homeomorphism onto $V(\mathfrak I)$.
\end{proof}

The next example shows that the requirement that $f_1 \circ g = f_2$ in Theorem~ \ref{thm: inj surj iso}\ref{thm: surj} cannot be omitted. Note that we do have $(f_1 \circ g)(y_2) \subseteq f_2(y_2)$ for all $y_2 \in  Y_2$, since $g$ defines a morphism of patch bundles, so what the example shows is that equality cannot be replaced with the weaker set inclusion criterion.

\begin{example}
Let $S$ be a proper subring of $R$.
We can view $S$ and $R$ as patch algebras for $R$ by taking $B_S=B_R=\{0,1\}$ and, under the canonical identification $R[\{0,1\}]\cong R$, taking $\iota_S$ to be the inclusion $S\hookrightarrow R$ and $\iota_R$ to be the identity on $R$.

The corresponding patch bundles are
\[
f_S\colon\; \{*\}\longrightarrow\Sigma(R): *\longmapsto S
\: \text{ and }\:
f_R\colon\; \{*\}\longrightarrow\Sigma(R): *\longmapsto R,
\]
where $\{*\}$ denotes the one-point space. The induced morphism of patch bundles $f_R\to f_S$ is given by the identity
on $\{*\}$. Its underlying map is a homeomorphism, but the morphism
$S\hookrightarrow R$ is not surjective.
\end{example}

\section[Patch ringed spaces]{Patch ringed spaces}\label{sec: patch ringed spaces}

Patch bundles and patch algebras give topological and algebraic realizations, respectively, of continuously indexed spaces of subrings.  
In this section, we round out this picture with a geometric realization given  by the notion of a patch ringed space.
To do so, we pass from patch bundles to ringed Stone spaces, see them as sheafifications of presheaves naturally arising from patch presheaves, and recover patch algebras as their rings of global sections. The main algebra-geometry duality is given in Theorem \ref{ASR}.

\subsection[Patch bundles and ringed spaces]{Patch bundles and ringed spaces}\label{subsec: patch bundles and ringed spaces} We denote by $\underline R_Y$ the constant sheaf on a Stone space $Y$ induced by a ring $R$, whose ring of sections on an open subset $U$ of $Y$ is   $$\underline R_Y(U)=C(U,R_{\mathrm{disc}}),$$
where $C(U,R_{\mathrm{disc}})$ denotes the ring of continuous functions from $U$ to $R$, where $R$ has the discrete topology.
\begin{definition} \label{patch ringed space def} \label{morphism for PRingSp} 
Let $R$ be a ring. A {\it patch ringed space for $R$} is a pair $(Y,\mathcal O_Y)$ consisting of a Stone space $Y$ and a subsheaf of rings $\mathcal O_Y\subseteq \underline R_Y,$ where $\underline R_Y$ denotes the constant sheaf on $Y$ associated to $R$.
{Denote by $\mathsf{PRS}_R$ the category whose objects are patch ringed spaces over $R$. A {\it morphism}
$(g,g^\#)\colon\; (Y_1,\mathcal O_{Y_1})\to(Y_2,\mathcal O_{Y_2})$ in $\mathsf{PRS}_R$ is a morphism of ringed spaces for which the following diagram commutes:
$$
\begin{tikzcd}
\mathcal O_{Y_2} \arrow[r,hook] \arrow[d,"g^\#"'] & \underline R_{Y_2} \arrow[d,"\overline{g^\#}"] \\
g_*\mathcal O_{Y_1} \arrow[r,hook] & g_*\underline R_{Y_1},
\end{tikzcd}
$$
where $\overline{g^\#}$ is the morphism whose map on sections over an open subset $U\subseteq Y_2$ sends $\sigma\in C(U,R_{\mathrm{disc}})$ to $\sigma\circ g\in C(g^{-1}(U),R_{\mathrm{disc}})$.}
\end{definition}

The stalks of $\underline R_Y$ are canonically isomorphic
to $R$ \cite[\href{https://stacks.math.columbia.edu/tag/007B}{Tag 007B}]{stacks-project}, so if $(Y,\mathcal O_Y)$ is a patch ringed space and $y \in Y$, the inclusion of sheaves $\mathcal{O}_Y\hookrightarrow \underline R_Y$ induces an injective map on stalks $\mathcal{O}_{Y,y}\hookrightarrow (\underline R_Y)_y$ (since a filtered colimit of injections is an injection) which, composed with the stalk identification, gives an embedding of $\mathcal O_{Y,y}$ into $R$. We denote the image of this map  by $R_y$ when it is clear from context which subsheaf $\O_Y$ we are working with.

\begin{lemma} \label{lem: rep}
Let $R$ be a ring, and let $(Y,\mathcal O_Y)$ be a patch ringed space for $R$. For every open subset $U\subseteq Y$,
\begin{eqnarray}\label{section eq}\mathcal O_Y(U)=
\left\{\sigma\in C(U,R_{\mathrm{disc}}):
\sigma(u)\in R_u \text{ for all } u\in U \right\}.\end{eqnarray}
\end{lemma}

\begin{proof}
Since $\mathcal O_Y$ is a subsheaf of $\underline R_Y$, $\mathcal O_Y(U) \subseteq C(U,R_{\mathrm{disc}})$.
 Let
\(\sigma\in\mathcal O_Y(U)\), and let \(u\in U\). The germ
$[U,\sigma]$ belongs to the stalk $\mathcal O_{Y,u}$. Under the canonical embedding $\mathcal O_{Y,u}\hookrightarrow R$, the germ $[U,\sigma]$ maps to $\sigma(u)$. Hence $\sigma(u)\in R_u.$ Therefore, the inclusion $\subseteq$ in  (\ref{section eq}) holds. 

Conversely, suppose that $\sigma\in C(U,R_{\mathrm{disc}})$
satisfies $\sigma(u)\in R_u$ for all $u\in U.$
We show that $\sigma\in\mathcal O_Y(U)$. Fix $u\in U$. Since
$\sigma(u)\in R_u$, there exists a germ $
\tau_u \in\mathcal O_{Y,u}$
whose image in $R$ is $\sigma(u)$. Choose an open neighborhood
$V\subseteq U$ of $u$ and a section
$
\tau \in\mathcal O_Y(V)$
representing this germ. Viewing $\tau$ as a section of
$\underline R_Y(V)=C(V,R_{\mathrm{disc}})$, we have $\tau(u)=\sigma(u).$
Since both $\tau$ and $\sigma|_V$ are in $C(V,R_{\mathrm{disc}})$, after replacing $V$ by a smaller open neighborhood of
$u$, we may assume  that both $\tau$ and $\sigma|_V$ are constant
on $V$ with constant value $\sigma(u)$. Thus
$\tau=\sigma|_V$. This proves that there is an open cover
$\{V_i\}_{i\in I}$ of $U$ such that $\sigma|_{V_i}\in \mathcal O_Y(V_i)$
for all $i$. On intersections $V_i\cap V_j$, these local sections agree because they are both equal to the restriction of the same function $\sigma$. Since $\mathcal O_Y$ is a sheaf, the local sections glue to a unique section of $\mathcal O_Y(U)$, which must then be $\sigma$. Hence $\sigma\in\mathcal O_Y(U),$
which proves the lemma.
\end{proof}

\begin{lemma}\label{stalk bundle} Let $R$ be a ring, and let $(Y,\mathcal O_Y)$ be a patch ringed space for $R$.
The map $$f\colon\; Y \longrightarrow \Sigma(R) : y\longmapsto R_y$$ is a patch bundle over $R$.
\end{lemma}
\begin{proof} { Fix $x\in R$, and let $U_x=\{R_y:y\in Y,\ x\in R_y\}.$}
This is a basic open subset of the subspace $\{R_y:y \in Y\}$ of $\Sigma(R)$. We claim that $f^{-1}({ U_x}) = \{y \in Y:x \in R_y\}$  is open in $Y$.   
Let $y \in f^{-1}({ U_x})$, so that $x\in R_y$. There is  an open neighborhood $V$ of $y$
and a section $\sigma$ in $\O_Y(V)$ such that $\sigma(y) =x$. 
Since $\sigma \in C(V,R_{\mathrm{disc}})$, there is an open neighborhood $W \subseteq V$ of $y$ such that $\sigma|_W$ has constant value $\sigma(y) =x$.
Since $\sigma|_W\in\mathcal O_Y(W)$, Lemma~\ref{lem: rep} implies that $x\in R_w$ for every $w \in W$. Hence
$W\subseteq f^{-1}({ U_x}).$
Thus every point of $f^{-1}({ U_x})$ has an open neighborhood contained in
$f^{-1}({ U_x})$, so $f^{-1}({ U_x})$ is open in $Y$ and $f$ is continuous. 
\end{proof}

\begin{definition} Let $R$ be a ring, and let $(Y,\mathcal O_Y)$ be a patch ringed space for $R$.
The {\it patch bundle of $(Y,\mathcal O_Y)$} is the continuous map $f\colon\; Y \rightarrow \Sigma(R) : y\longmapsto R_y$ from Lemma~\ref{stalk bundle}, and the {\it patch presheaf of $(Y,\mathcal O_Y)$} is the patch presheaf associated to this patch bundle. 
\end{definition}

\begin{definition}\label{def: bundle ringed functors}
We define the functors
\[
\mathcal S\colon\; \mathsf{PBun}_R\longrightarrow\mathsf{PRS}_R
\quad\text{and}\quad
\mathcal S^{-1}\colon\; \mathsf{PRS}_R\longrightarrow\mathsf{PBun}_R
\]
as follows.
On objects, $\mathcal S$ sends a patch bundle $f\colon\; Y\to\Sigma(R)$ to $\mathcal S(f)=(Y,\mathcal O_f)$, where
$$\mathcal O_f(U)=\{\sigma\in C(U,R_{\mathrm{disc}}):
\sigma(y)\in f(y)\text{ for every }y\in U\},
$$
for every open subset $U\subseteq Y$, with the usual restriction maps.
On morphisms, $\mathcal S$ sends a bundle morphism $f_1\to f_2$ represented by $g\colon\; Y_1\to Y_2$ to
\[
\mathcal S(g)=(g,g^\#)\colon\; (Y_1,\mathcal O_{f_1})\longrightarrow(Y_2,\mathcal O_{f_2}),
\]
where, for every open subset $U\subseteq Y_2$, 
\[
g^\#_U\colon\; \mathcal O_{f_2}(U)\longrightarrow\mathcal O_{f_1}(g^{-1}(U))\colon\; \sigma\longmapsto\sigma\circ g.
\]
On objects, $\mathcal S^{-1}$ sends a patch ringed space $(Y,\mathcal O_Y)$ to its stalk bundle
\[
\mathcal S^{-1}(Y,\mathcal O_Y)\colon\;
Y\longrightarrow\Sigma(R)\colon\; y\longmapsto R_y.
\]
On morphisms, it sends
$(g,g^\#)\colon\; (Y_1,\mathcal O_{Y_1})\to(Y_2,\mathcal O_{Y_2})$
to the bundle morphism
\[
\mathcal S^{-1}(g,g^\#)\colon\; \mathcal S^{-1}(Y_1,\mathcal O_{Y_1})
\longrightarrow\mathcal S^{-1}(Y_2,\mathcal O_{Y_2})
\]
represented by the same continuous map $g\colon\; Y_1\to Y_2$.
\end{definition}

\begin{theorem}\label{thm: bundle ringed equivalence}
The functors $\mathcal S$ and $\mathcal S^{-1}$ are quasi-inverse equivalences between $\mathsf{PBun}_R$ and $\mathsf{PRS}_R$.
\end{theorem}
\begin{proof}
Let $f:Y \rightarrow \Sigma(R)$ be a patch bundle.  It is clear that $\mathcal O_f$ is a subsheaf of $\underline R_Y$. We claim that the stalk at $y$, viewed as a subring of $R$, is $f(y)$. To see this is the case, observe first that evaluation at $y$ sends every germ in $(\mathcal O_f)_y$
to an element of $f(y)$. Conversely, if $r\in f(y)$, the constant function with value $r$ on the open neighborhood $f^{-1}(\mathcal U(r))$ is a section of $\mathcal O_f$ representing $r$ at $y$.
If $g$ represents a morphism $f_1\to f_2$ of patch bundles, then
$$
(\sigma\circ g)(y)\in f_2(g(y))\subseteq f_1(y)
$$
for every section $\sigma$ of $\mathcal O_{f_2}$. Thus composition with $g$ defines the required morphism of patch ringed spaces. Conversely, the commutative diagram in Definition~\ref{morphism for PRingSp}, after passing to stalks and evaluating in $R$, gives the inclusions
$
R_{2,g(y)}\subseteq R_{1,y}.
$
Therefore, by Lemma~\ref{stalk bundle}, the assignments defining $\mathcal S^{-1}$  give objects and morphisms of $\mathsf{PBun}_R$. 

Both assignments keep the underlying continuous maps, so
they preserve identities and compositions on the underlying spaces. On sections, the same follows from $\sigma\circ\operatorname{id}=\sigma$ and $(\sigma\circ h)\circ g=\sigma\circ(h\circ g)$ for composable continuous maps $g$ and $h$.

For every $y\in Y$ , the image of the stalk $(\mathcal{O}_f)_y$ under evaluation at $y$ is $f(y)$. Thus the patch bundle associated with $(Y,\mathcal O_f)$ is $f$ , so $\mathcal S^{-1}\mathcal S(f)=f$. This identification is natural, since both functors preserve the underlying maps of patch bundles.

On the other hand, let $(Y,\mathcal O_Y)$ be a patch ringed space, and let $f_{\mathcal O}$ be its associated patch bundle, given by
$f_{\mathcal O}(y)=R_y$. By Lemma~\ref{lem: rep}, for every open
subset $U\subseteq Y$,
\[
\mathcal O_{f_{\mathcal O}}(U)
=
\{\sigma\in C(U,R_{\mathrm{disc}}):
\sigma(u)\in R_u\text{ for all }u\in U\}
=\mathcal O_Y(U).
\]
Thus $\mathcal S\mathcal S^{-1}(Y,\mathcal O_Y)
=(Y,\mathcal O_Y)$. This identification is also natural, since the maps on sections are determined by
precomposition with the underlying continuous maps. Hence both composites are naturally isomorphic to the corresponding identity
functors.
\end{proof}

\subsection[Sheafification and global sections]{Sheafification and global sections}\label{subsec: sheafification and global sections}

As observed in \cite[Section~3.2]{CosOlb}, the terminology ``patch presheaf'' reflects the fact that these objects can be viewed as ordinary presheaves on Stone spaces. The following lemma makes this precise.

\begin{lemma} \label{lem: presheaf from patch} Let $R$ be a ring, and let $({{\mathcal R}}=\{R_e\subseteq R:e\in B\}, B)$ be a patch presheaf for $R$. There is a presheaf $\O_{{\mathcal R}}$ on the Stone space $\Max(B)$ such that $\O_{{\mathcal R}}(U_e) = R_e$ for every $e\in B$. 
\end{lemma} 
\begin{proof} 
First we define $\O_{{\mathcal R}}$ on the clopen sets of   $\Max(B)$ by setting  $\O_{{\mathcal R}}(U_e) = R_e$ for each $e \in B$.
If $e,f \in B$ and $U_e \subseteq U_f$, then $e\le f$  by Stone duality and hence $$R_f = \O_{\mathcal R}(U_f) \subseteq \O_{\mathcal R}(U_e)=R_e.$$ The inclusion of $\O_{\mathcal R}(U_f)$ in $\O_{\mathcal R}(U_e)$ is then the restriction morphism.   The above construction can be extended from clopen sets to all open sets of $\Max(B)$ in a standard way in which the ring of sections $\O_{{\mathcal R}}(U)$ over an arbitrary open set $U$ is the inverse limit of the $R_e$ over all $e\in B$ with $U_e\subseteq U$.
\end{proof}

\begin{theorem}\label{thm: patch ringed space}
Let $R$ be a ring, and let $(Y,\mathcal O_Y)$ be a patch ringed space for $R$.
The sheaf $\O_Y$ is the sheafification of the presheaf $\O_{\mathcal R}$, where $ (\mathcal R,\Clop(Y))$ is the patch presheaf of $(Y,\mathcal O_Y)$.
\end{theorem}
\begin{proof} Let $f\colon\; Y \rightarrow \{R_y:y \in Y\}$ denote the patch bundle of $(Y,\mathcal O_Y)$, and consider the patch presheaf of the patch bundle $f$,
$$(\mathcal R,\Clop(Y))= {\mathcal P}(f) = (\{R_U:U \in {\rm Clop}(Y)\},{\rm Clop}(Y)),$$  where, as recalled in Subsection \ref{subsect: patch presheaves and bundles}, for each clopen set $U$ of $Y$, $R_U$ is the intersection of the rings in $f(U)$.
As observed above, $\O_{\mathcal R}$ is a presheaf on the Stone space $Y$, and by Stone duality, $Y$ is homeomorphic to the space of maximal ideals of the Boolean algebra ${\rm Clop}(Y)$. 
The presheaf $\O_{\mathcal R}$ is defined by $\O_{\mathcal R}(U) = R_U$ for all clopen sets $U$ of $Y$.  
Following \cite[\href{https://stacks.math.columbia.edu/tag/007X}{Tag 007X}]{stacks-project}, the sheafification $\O_{\mathcal R}^\#$ of $\O_{\mathcal R}$ is defined for each open subset $U$ of $Y$ by letting $\O_{\mathcal R}^\#(U)$ be the set of all functions 
$\sigma\colon\; U \rightarrow R$ such that
\begin{enumerate}[label=\textup{(\roman{*})}, mode=unboxed]
\item\label{enum a} $\sigma(u) \in f(u) = {{\mathcal R}}_u$ for all $u$ in $U$, and
\item\label{enum b} { for every $u\in U$, there exist a clopen neighborhood $u\in V\subseteq U$ and $t\in R_V$ such that $\sigma(v)=t$ for every $v\in V$.}  
\end{enumerate}

Let $U\subseteq Y$ be open. By Lemma~\ref{lem: rep}, the sections of $\mathcal O_Y(U)$ are  the locally constant functions with values in the rings $f(u)$. Therefore, conditions~\ref{enum a} and~\ref{enum b}  imply that $\mathcal O_{\mathcal R}^{\#}(U)\subseteq\mathcal O_Y(U)$. Conversely, let $\sigma\in\mathcal O_Y(U)$ and fix $u\in U$. Since $\sigma$ is locally constant and the clopen subsets form a basis, there is a clopen neighborhood $u\in V\subseteq U$ on which $\sigma$ has constant value $t=\sigma(u)$. For every $v\in V$, $t=\sigma(v)\in f(v)$, so $t\in R_V$. Thus $\sigma$ satisfies~\ref{enum b}, while~\ref{enum a} follows from Lemma~\ref{lem: rep}. Hence $\mathcal O_{\mathcal R}^{\#}(U)=\mathcal O_Y(U)$ for every open $U$, proving the theorem.
\end{proof}

In particular, the structure sheaf of \(\mathcal S(f)\) is the sheafification of the presheaf associated with the patch presheaf \(\mathcal P(f)\).
Given a patch presheaf $({\mathcal R}=\{R_e:e\in B\},B)$, let $Y=\Max(B)$ and let $\O_{{\mathcal R}}$ be the associated presheaf on $Y$. By Theorems~\ref{thm: bundle presheaf equivalence} and~\ref{thm: patch ringed space} and Definition~\ref{def: bundle ringed functors}, for every open subset $U\subseteq Y$ we have 
\begin{equation}\label{eq: sheafification}
    \mathcal O_{{\mathcal R}}^{\#}(U) = \left\{\sigma\in C(U,R_{\mathrm{disc}}): \sigma(\mathfrak m)\in{\mathcal R}_{\mathfrak m}\text{ for every }\mathfrak m\in U \right\}.
\end{equation}

\begin{corollary}\label{cor: sheafification equivalence}
Sheafification defines a category equivalence
$$ \mathsf{PPre}_R\longrightarrow\mathsf{PRS}_R: (\mathcal R,B)\longmapsto
\bigl(\Max(B),\mathcal O_{\mathcal R}^{\#}\bigr)$$
that is  naturally isomorphic to $\mathcal S\circ\mathcal P^{-1}$.
\end{corollary}
\begin{proof}
For each patch presheaf $(\mathcal R,B)$, the description of \(\mathcal O_{\mathcal R}^{\#}(U)\) given in \eqref{eq: sheafification} and Definition~\ref{def: bundle ringed functors} give a canonical isomorphism
$$
\bigl(\Max(B),\mathcal O_{\mathcal R}^{\#}\bigr)
\cong (\mathcal S\circ\mathcal P^{-1})(\mathcal R,B)
$$
that  is the identity on the underlying space and, on sections, sends each $R$-valued function to itself. 
We claim that this isomorphism is also compatible with morphisms.
Let $\psi\colon\; (\mathcal R_1,B_1)\to(\mathcal R_2,B_2)$ be a morphism of patch presheaves, where $\mathcal R_i=\{R_{i,e}:e\in B_i\}$. Let $h\colon\; B_2\to B_1$ be the Boolean homomorphism defining $\psi$, and let $g\colon\; Y_1\to Y_2$ be its Stone-dual map, where $Y_i=\Max(B_i)$. For each $e\in B_2$, the inclusion $R_{2,e}\subseteq R_{1,h(e)}$ sends a constant section with value $r$ on $U_e$ to the constant section with the same value on $g^{-1}(U_e)=U_{h(e)}$. Thus the induced map after sheafification sends a section $\sigma$ over an open subset $U\subseteq Y_2$ to $\sigma\circ g$ over $g^{-1}(U)$, which is also the map in Definition~\ref{def: bundle ringed functors}.
Therefore, with  $f_i=\mathcal P^{-1}(\mathcal R_i,B_i)$, we  have, for every open subset $U\subseteq Y_2$, the commutative diagram
$$
\begin{tikzcd}[column sep=large, row sep=large]
\mathcal O_{\mathcal R_2}^{\#}(U)
  \arrow[r, "\sim"] \arrow[d, "{\sigma\mapsto\sigma\circ g}"']
& \mathcal O_{f_2}(U)
  \arrow[d, "{\sigma\mapsto\sigma\circ g}"] \\
\mathcal O_{\mathcal R_1}^{\#}(g^{-1}(U))
  \arrow[r, "\sim"]
& \mathcal O_{f_1}(g^{-1}(U))
\end{tikzcd}
$$
where each horizontal arrow sends an $R$-valued function to itself.
Consequently, the canonical isomorphisms of patch ringed spaces $(Y,\mathcal O_{\mathcal R}^{\#})\cong(Y,\mathcal O_f)$  commute with morphisms and define a natural isomorphism between the sheafification functor and $\mathcal S\circ\mathcal P^{-1}$. Since $\mathcal P^{-1}$ and $\mathcal S$ are equivalences by Theorems~\ref{thm: bundle presheaf equivalence} and~\ref{thm: bundle ringed equivalence}, the  sheafification functor is also an equivalence.
\end{proof}

We next construct an isomorphism from the  patch algebra of a patch bundle to  the global sections of its patch ringed space.
\begin{definition}\label{lem: global section algebra}
Let $R$ be a ring, and let $(Y,\mathcal O_Y)$ be a patch ringed space for $R$.
The {\it patch algebra} of $(Y,\mathcal O_Y)$ is the ring $A_Y=\O_Y(Y)$ of global sections of $\O_Y$  equipped with the Specker presentation
$
B_Y=\{1_U:U\in\Clop(Y)\} $ and $
\iota_{A_Y}\colon\; A_Y\hookrightarrow R[B_Y],
$
where $\iota_{A_Y}$ is the restriction of the canonical isomorphism $C(Y,R_{\mathrm{disc}})\to R[B_Y]$ fixing $B_Y$.
\end{definition}

That this defines a patch algebra follows from Lemma~\ref{lem: rep} and Theorem~\ref{C rep}.

\begin{theorem}\label{thm: patch ringed space 2}
Let $R$ be a ring, let $(Y,\mathcal O_Y)$ be a patch ringed space for $R$, and let $f$ be its associated patch bundle. For every clopen subset $U\subseteq Y$, the  ring of sections $\mathcal O_Y(U)$, equipped with its canonical Specker presentation, is canonically isomorphic to the patch algebra $A(f|_U)$.
Moreover, with $A=\mathcal O_Y(Y)$ and $e_U=1_U$ the identity of $e_UA$, there are canonical ring isomorphisms
$$
\mathcal O_Y(U)\cong e_UA\cong A/(1-e_U)A.
$$
\end{theorem}
\begin{proof}
We first prove the assertion for $U=Y$. By Lemma~\ref{lem: rep},   $$\mathcal O_Y(Y) =
\left\{\sigma \in C(Y,R_{\mathrm{disc}}): \sigma(y)\in R_y \text{ for all } y\in Y \right\}.$$
Since $R_y=f(y)$, Theorem~\ref{C rep} identifies this ring with $A(f)$ as a patch algebra, with the Specker presentation of Definition~\ref{lem: global section algebra}.

For a clopen subset $U\subseteq Y$, the pair $(U,\mathcal O_Y|_U)$ is a patch ringed space with associated bundle $f|_U$. Applying the case just proved yields the canonical isomorphism $\mathcal O_Y(U)\cong A(f|_U)$ of patch algebras, with the Specker presentation in  Definition~\ref{lem: global section algebra} applied to $(U,\mathcal O_Y|_U)$.
Finally, every section over $U$ extends   to a global section that is $0$ off $U$. These extensions form $e_UA$, and the kernel of restriction $A\to\mathcal O_Y(U)$ is $(1-e_U)A$.
\end{proof}

The global-sections functor $\Gamma\colon\;\mathsf{PRS}_R\to\mathsf{PAlg}_R^{\mathrm{op}}$ is defined in the usual way. A patch ringed space $(Y,\O_Y)$ is mapped to its patch algebra of global sections. A morphism
\[
(g,g^\#)\colon\; (Y_1,\O_{Y_1})\longrightarrow(Y_2,\O_{Y_2})
\]
is mapped to the opposite of the patch algebra morphism
\[
g^\#_{Y_2}\colon\; \O_{Y_2}(Y_2)\longrightarrow\O_{Y_1}(Y_1).
\]
That this is a morphism of patch algebras follows from Definition~\ref{morphism for PRingSp} and the Specker presentation given in Definition~\ref{lem: global section algebra}.

\begin{theorem} \label{ASR}
The global-sections functor 
$$\Gamma\colon\; \mathsf{PRS}_R\longrightarrow \mathsf{PAlg}_R^{\mathrm{op}}\colon\; \:(Y,\mathcal O_Y)\longmapsto A_Y$$
with the Specker presentation of  Definition~\ref{lem: global section algebra} is naturally isomorphic to $\mathcal A\circ\mathcal S^{-1}$ and is an equivalence with quasi-inverse $\mathcal S\circ\mathcal A^{-1}$. Thus $\mathsf{PRS}_R$ is equivalent to each of $\mathsf{PAlg}_R^{\mathrm{op}}$, $\mathsf{PPre}_R$, and $\mathsf{PBun}_R$.
\end{theorem}
\begin{proof}
Theorem~\ref{thm: patch ringed space 2} gives an isomorphism of  $(\mathcal A\circ\mathcal S^{-1})(Y,\mathcal O_Y)$ with the patch algebra of global sections. This isomorphism is natural because, for a morphism $(g,g^\#)$ in $\mathsf{PRS}_R$, Theorem~\ref{C rep} implies that both the induced algebra morphisms and the maps on global sections send $\sigma$ to $\sigma\circ g$.
Therefore, $\Gamma$ is naturally isomorphic to $\mathcal A\circ\mathcal S^{-1}$. By Theorems~\ref{thm: bundle ringed equivalence} and~\ref{cor: PAlg PBundle}, $\mathcal A\circ\mathcal S^{-1}$ is a category equivalence with quasi-inverse $\mathcal S\circ\mathcal A^{-1}$. Since $\Gamma$ is naturally isomorphic to $\mathcal A\circ\mathcal S^{-1}$, it is also an equivalence, with quasi-inverse $\mathcal S\circ\mathcal A^{-1}$.
\end{proof}

Together with the bundle-presheaf equivalence of Theorem~\ref{thm: bundle presheaf equivalence}, Theorem~\ref{ASR} completes the following diagram of category equivalences, which commutes up to natural isomorphism.
$$
\begin{tikzcd}[column sep=large,row sep=large]
\mathsf{PBun}_R
 \arrow[r,"\mathcal P"]
 \arrow[d,"\mathcal S"']
 \arrow[dr,"\mathcal A"{description}]
&
\mathsf{PPre}_R
 \arrow[d,"\mathcal A\circ\mathcal P^{-1}"]
\\
\mathsf{PRS}_R
 \arrow[r,"\Gamma"']
&
\mathsf{PAlg}_R^{\mathrm{op}}
\end{tikzcd}
$$

 Let $f\colon\; Y\to\Sigma(R)$ be a patch bundle, and let $\mathcal O_f$ be its associated sheaf. By Theorem~\ref{thm: patch ringed space 2}, the ring of sections of $\O_f$ on a clopen subset $U$ is canonically isomorphic to the patch algebra $A(f|_U)$. For nonempty $U$, 
 the presheaf coefficient ring $R_{U}$ embeds in this algebra as the constant sections. The following examples illustrate the relationship between these two rings.

\begin{example}
Let $X$ be a finite collection of subrings of $R$, and consider the tautological patch bundle $f$ on $X^{\rm patch}$ from Example~\ref{finite 1}. For every nonempty subset $E\subseteq X$, the presheaf coefficient ring and the ring of sheaf sections are
$$ R_E=\bigcap_{T\in E}T \: \:\text{ and } \: \mathcal O_f(E)\cong A(f|_E)\cong\prod_{T\in E}T. $$
The canonical map from $R_E$ to this product is the diagonal map $r\mapsto(r)_{T\in E}$.  In particular, if $X=\{S,T\}$ with $S\subsetneq T$, then
$ R_X=S$ and $\mathcal O_f(X)\cong S\times T.$

More generally, consider the patch bundle $f\colon\; Y\to\{S,T\}$ of Example~\ref{finite 1}, and let $U=f^{-1}(T)$. The subset $U$ is open but need not be clopen. For every nonempty clopen $E\subseteq Y$, $ R_E=T$ if $E \subseteq U$ and $R_E=S$ if $E \not \subseteq U$. Thus $R_E$ identifies with the ring of
constant sections on $E$, whereas
$\mathcal O_f(E)=
\{\sigma\in C(E,T_{\mathrm{disc}}):
\sigma(E\setminus U)\subseteq S\}.$
In particular, sheaf sections may take different values on different clopen subsets of $E$. 
\end{example}

\begin{example}
Consider the patch bundle $f\colon\; Y=\mathbb N\cup\{\infty\}\to\Sigma(\mathbb Q)$ of Example~\ref{chain 1}, where $f(n)=S_n$ and $f(\infty)=\mathbb Q$. For the clopen  $ E_m=\{m,m+1,\ldots\}\cup\{\infty\}, $ the presheaf coefficient ring is $R_{E_m}=S_m$, whereas
$$ \mathcal O_f(E_m)\cong A(f|_{E_m})\cong
\left\{
(a_n)_{n\ge m}\in\prod_{n\ge m}S_n: (a_n)\text{ is eventually constant} \right\}. $$
The canonical embedding of $S_m$ into $\mathcal O_f(E_m)$ sends each element to the corresponding constant sequence. If instead $E\subseteq\mathbb N$ is nonempty and finite, then $ R_E=S_{\min E} $ and 
$\mathcal O_f(E)\cong\prod_{n\in E}S_n. $
Every clopen subset containing $\infty$ is the disjoint union of a tail and finitely many isolated points, so these descriptions determine the section rings on all nonempty clopens.
\end{example}

\begin{example}
Consider the tautological patch bundle $\iota\colon\; \operatorname{Zar}(K)^{\rm patch}\to\Sigma(K)$ of Example~\ref{ex: zariski inclusion}. For a nonempty clopen subset $E\subseteq\operatorname{Zar}(K)^{\rm patch}$, the presheaf coefficient ring is $ R_E=\bigcap_{W\in E}W, $ whereas the  ring of sections of $E$ for the sheaf is $$ \mathcal O_\iota(E) \cong  A(\iota|_E) \cong 
\{\sigma\in C(E,K_{\mathrm{disc}}):   \sigma(W)\in W\text{ for every }W\in E\}. $$
An element of $R_E$ must belong to every valuation ring in $E$. A sheaf section can instead take different values on different clopen pieces of $E$,
provided each value belongs to every valuation ring in the corresponding piece.

Now let $F$ be a field extension of $K$, and consider the reindexed patch bundle $\iota\circ\phi$, where $\phi(V)=V\cap K$. For a nonempty clopen subset
$E\subseteq\operatorname{Zar}(F)^{\rm patch}$, the comparison becomes $R_E=\bigcap_{V\in E}(V\cap K),$ and
$$ \mathcal O_{\iota\circ\phi}(E)= \{\sigma\in C(E,K_{\mathrm{disc}}):   \sigma(V)\in V\cap K\text{ for every }V\in E\}.$$
If $E$ contains distinct valuation rings $V_1,V_2$ with the same contraction to $K$, choose a clopen subset $D\subseteq E$ containing $V_1$ but not $V_2$. Then $1_D$ is a sheaf section taking different values at these points, although their stalks are the same subring of $K$.
Thus sheaf sections can distinguish distinct valuation rings with the same contraction to $K$, even though their stalks coincide.
Constant sections, by contrast, cannot distinguish these points.
\end{example}

We postpone the discussion of the example arising from integral separated quasi-compact schemes (Example~\ref{scheme 1}) to Appendix~\ref{sec: w localization}, where we relate the associated patch algebra to Bhatt-Scholze $w$-localization.

\appendix

\section[Pierce duality]{Pierce duality}\label{sec: Pierce} In the case in which the ring $R$ is indecomposable, the category equivalence between $\mathsf{PAlg}_R^{\mathrm{op}}$ and $\mathsf{PRS}_R$  can be viewed as an instance of Pierce duality. We explain how to do this in this appendix. 

First we recall Pierce duality for commutative rings from Pierce's monograph \cite{Pierce}. 
The  {\it  Pierce sheaf} of a commutative ring $A$ is the
sheafification $ \Pi_A$ of the presheaf on $\Max(\Id(A))$ defined on clopens $U_e$ by $U_e\mapsto eA$ for $e\in\Id(A)$,
and whose restriction map $eA\to fA$, for $f\le e$, is multiplication
by $f$. The stalks of $\Pi_A$ are called the {\it Pierce stalks} of $A$. For
$\mathfrak m\in\Max(\Id(A))$, the stalk at $\mathfrak m$
is canonically isomorphic to $A/\mathfrak mA$. (For more details on this, see \cite[pp. 233-234]{BV}.)

Let $\mathsf{CRing}$ denote the category of commutative rings and unital ring homomorphisms, and let $\mathsf{Pierce}$ denote the category whose objects are ringed Stone spaces with nonzero indecomposable stalks and whose morphisms are morphisms of ringed spaces. Pierce duality \cite[Theorem~6.6, p.~24]{Pierce} is the category equivalence
$$
\mathcal D\colon\;\mathsf{CRing}^{\mathrm{op}}\longrightarrow\mathsf{Pierce}: A\longmapsto\bigl(\Max(\Id(A)),\Pi_A\bigr).
$$

For indecomposable $R$, the next theorem shows that our duality between patch algebras and patch ringed spaces agrees, up to natural isomorphism, with classical Pierce duality after forgetting the Specker presentations and the specified inclusions into constant sheaves.
 
\begin{theorem}\label{Pierce duality restriction}
Assume $R$ is indecomposable. Let $(Y,\mathcal O_Y)$ be a patch ringed space for $R$, and let $A=\mathcal O_Y(Y)$ be its patch algebra. Under the canonical homeomorphism
$$
Y\cong\Max(\operatorname{Id}(A)),
$$
the structure sheaf $\mathcal O_Y$ is canonically isomorphic to the Pierce sheaf $\Pi_A$. The resulting isomorphisms of ringed spaces
$$
(Y,\mathcal O_Y)\cong
\bigl(\Max(\operatorname{Id}(A)),\Pi_A\bigr)
$$
are natural with respect to morphisms of patch ringed spaces.
\end{theorem}

\begin{proof} Let $B_Y=\{1_U:U\in\Clop(Y)\}$, and let $A=\O_Y(Y)$ be the patch algebra of $(Y,\O_Y)$ 
with the embedding $\iota_A\colon\; A\hookrightarrow R[B_Y]$ from  Definition~\ref{lem: global section algebra}.
Since $R$ is indecomposable, the idempotents of
$C(Y,R_{\mathrm{disc}})$ are exactly the characteristic functions of clopen subsets of $Y$. Hence $\operatorname{Id}(A)=B_Y.$ Thus the Pierce space of $A$ is $\Max(B_Y)$, which by Stone duality is canonically
homeomorphic to $Y$.
Under this homeomorphism, Theorem~\ref{thm: patch ringed space 2} gives an isomorphism of $\O_Y(U_e)$ with $eA$ via the map that extends each section by zero outside $U_e$, for every $e\in B_Y$.
Here $e$ is the characteristic function of $U_e$, so $eA$ consists precisely of the global sections that vanish outside $U_e$. 
For $f\le e$, restriction to $U_f$ corresponds to multiplication by $f$.
Since the isomorphisms $\mathcal O_Y(U_e)\cong eA$ commute with restrictions, they induce an isomorphism of  $\mathcal O_Y$ with the Pierce sheaf $\Pi_A$. Hence $(Y,\O_Y)$ is canonically isomorphic to the Pierce ringed space of its patch algebra.

We now prove naturality. Let
$$
(g,g^\#)\colon\; (Y_1,\mathcal O_{Y_1})\longrightarrow (Y_2,\mathcal O_{Y_2})
$$
be a morphism of patch ringed spaces, and write $A_i=\mathcal O_{Y_i}(Y_i)$ and $\psi=g^\#_{Y_2}\colon\; A_2\to A_1$. By Definition~\ref{morphism for PRingSp}, $g^\#(\sigma)
 =\sigma\circ g$.

 Thus, for $e\in\operatorname{Id}(A_2)$, $\psi(e)$ is the characteristic function of $g^{-1}(U_e)$, so the Stone-dual map induced by $\psi$ corresponds to $g$ under the canonical homeomorphisms.
For $\sigma\in\mathcal O_{Y_2}(U_e)$, let $\widetilde{\sigma}$ be the section on $Y_2$ that equals $\sigma$ on $U_e$ and is zero elsewhere. Since $g^\#_{U_e}(\sigma)=\sigma\circ g$ and $\psi(\widetilde{\sigma})=\widetilde{\sigma}\circ g$, we have
$$ g^\#_{U_e}(\sigma) =\left.\psi(\widetilde{\sigma})\right|_{g^{-1}(U_e)}.$$
Since $\widetilde{\sigma}=e\widetilde{\sigma}$, we have $\psi(\widetilde{\sigma})\in\psi(e)A_1$. Thus the following diagram commutes, where the horizontal arrows are the extension-by-zero isomorphisms:
$$
\begin{tikzcd}[column sep=large, row sep=large]
\mathcal O_{Y_2}(U_e)
  \arrow[r, "\sim"] \arrow[d, "g^\#_{U_e}"']
& eA_2 \arrow[d, "{\psi|_{eA_2}}"] \\
\mathcal O_{Y_1}(g^{-1}(U_e))
  \arrow[r, "\sim"]
& \psi(e)A_1
\end{tikzcd}
$$
The restrictions $\psi|_{eA_2}\colon\; eA_2\to\psi(e)A_1$ are the maps that define the morphism of Pierce presheaves, so after sheafification 
the sheaf isomorphisms $\mathcal O_{Y_i}\cong\Pi_{A_i}$ constructed in the first part of the proof, for $i=1,2$, commute with the morphism induced by $\psi$. This proves naturality.
\end{proof}

For more on the connection between Pierce duality and (locally) Specker algebras, see \cite{BMO}.

\section[Bhatt-Scholze w-localization]{Bhatt-Scholze $w$-localization}\label{sec: w localization}

In \cite{BhattScholze}, Bhatt and Scholze introduce $w$-localization as part of their study of the pro-\'etale topology. For a quasi-compact quasiseparated scheme $X$, they construct   an affine $w$-local scheme $X^{\rm wl}$ that maps surjectively onto $X$. In this appendix we show that when $X$ is also separated and integral, we can realize this affine scheme as the spectrum of the patch algebra of the tautological patch bundle constructed in Examples~\ref{scheme 1} and~\ref{scheme 2}.

Following \cite{BhattScholze}, a spectral space $X$ is {\it weakly local ($w$-local)} if its set of closed points is closed and every connected component contains a unique closed point.

\begin{lemma}\label{semiregular} The following are equivalent for a ring $A$ with Jacobson radical $J=J(A)$. 
\begin{enumerate}[label=\textup{(\arabic{*})}, mode=unboxed]
\item\label{semiregular 1} $\Spec(A)$ is $w$-local.

\item\label{semiregular 2} $\dim(A/J) = 0$ and each prime ideal of $A$ is contained in a unique maximal ideal.  

\item\label{semiregular 3} $A/J$ is von Neumann regular and idempotents of $A/J$ lift to $A$. 
\end{enumerate}
\end{lemma}

\begin{proof}
\ref{semiregular 1} $\Rightarrow$ \ref{semiregular 2}: This follows from \cite[Lemma~2.1.4]{BhattScholze}.

\ref{semiregular 2} $\Rightarrow$ \ref{semiregular 3}: Since $\dim(A/J) =0$ and $A/J$ is reduced, $A/J$ is von Neumann regular. By \cite[Theorem~1.2]{DeMarcoOrsatti}, there is a continuous retraction $\phi\colon\; \Spec(A)\rightarrow \Max(A)=\Spec(A/J).$
For an idempotent $\bar e\in A/J$, the clopen subset $\phi^{-1}(U(\bar e))$,  corresponds to an idempotent $e\in A$. Since $\phi$ restricts to the identity on $\Spec(A/J)$, we have $U(e+J)=U(\bar e)$, and hence $e+J=\bar e$. Thus idempotents lift modulo $J$.

\ref{semiregular 3} $\Rightarrow$ \ref{semiregular 1}: Since $A/J$ is von Neumann regular, $\Max(A)=V(J)$ is closed in $\Spec(R)$. Since distinct maximal ideals are separated by idempotents of $A/J$,  lifting these idempotents separates distinct maximal ideals of $A$ by clopen subsets of $\Spec(A)$, so each connected component contains at most one maximal ideal since a clopen subset separating two maximal ideals would disconnect any connected component containing both.
which lies in the same connected component since $V(\mathfrak p)$ is connected for every prime
ideal $\mathfrak p$. Thus $\Spec(A)$ is $w$-local.
\end{proof} 

The rings satisfying the equivalent conditions of the lemma are precisely
the {\it semiregular rings} introduced by Nicholson \cite[Theorem~2.9]{Nicholson1976}.
We show next that many naturally occurring patch algebras are semiregular rings. Recall from subsection \ref{subsec: patch bundles}, that a patch bundle $f\colon\; Y \rightarrow \Sigma(R)$ is a spectral map (i.e., preimages of quasi-compact open sets are quasi-compact) if and only if for each $r \in R$, the set $\{y \in Y:r \in f(y)\}$ is a clopen of $Y$.

\begin{proposition} Let  $A$ be the patch algebra of a patch bundle $f\colon\; Y \rightarrow \Sigma(R)$. If $f$ is a spectral map whose image consists of local subrings of $R$, then $A$ is a semiregular ring and hence $\Spec(A)$ is $w$-local. 
\end{proposition}

\begin{proof} We represent $A$ as in Theorem~\ref{C rep}. If $e\in A$ is idempotent, then for every $y\in Y$, the element $e(y)$ is an idempotent of $f(y)$. Since $f(y)$ is local, we have $e(y)=0$ or $e(y)=1$. Since $e$ is continuous, $e^{-1}(1)$ is clopen, and it  follows that the idempotents of $A$ are the characteristic functions of clopen subsets of $Y$. Thus, by Proposition~\ref{prop: stalks as factor rings}, the Pierce stalks of $A$ are canonically isomorphic to the local rings $f(y)$, for $y\in Y$. Since all its Pierce stalks are local, $A$ is clean by \cite[Proposition~1.2]{BurSte79}. Therefore, for every prime ideal $\mathfrak p$ of $A$, the quotient $A/\mathfrak p$ is a clean domain, hence a local ring, since each of its elements $x$ has the property that either $x$ or $1-x$ is a unit.
Therefore, every prime ideal of $A$ is contained in a unique maximal ideal.

Let $J$ be the Jacobson radical of $A$. By Lemma~\ref{semiregular}, it remains to show
that $\dim(A/J)=0$. We do so by proving that $A/J$ is von Neumann regular.
For this, it suffices to show that the image of every $\sigma\in A$ in $A/J$ is a product
of an idempotent and a unit.
Fix $\sigma\in A$. We claim that
$$E=\{y\in Y:\sigma(y)\text{ is a unit in }f(y)\}$$
is clopen in $Y$. Choose a finite clopen partition $\{U_1,\ldots,U_n\}$ of $Y$ such that $\sigma$ is constant on each $U_i$, with value $r_i$. If $r_i$ is not a unit in $R$, then
$E\cap U_i=\varnothing$, since $f(y)$ is a subring of $R$ for every $y\in Y$. Otherwise, since $r_i\in f(y)$ for every $y\in U_i$, we have $E\cap U_i =U_i\cap f^{-1}(\mathcal U(r_i^{-1}))$, which is clopen because $f$ is spectral. Thus $E\cap U_i$ is clopen in either case, and since the partition is finite,
so is $E$.

Let $e=1_E$ and $u=e\sigma+(1-e)$. To complete the proof, it suffices to show that that $u$ is a unit in $A$ and $\sigma-eu\in J$. 
For every $y\in Y$, we have
\[
u(y)=
\begin{cases}
\sigma(y) & \text{if }y\in E,\\
1 & \text{if }y\notin E.
\end{cases}
\]
Thus $u(y)$ is a unit in $f(y)$ for every $y\in Y$. Since $u$ is locally constant, so is the function $y\mapsto u(y)^{-1}$, which belongs to $A$ and is the inverse of $u$.

Moreover, since $eu=e\sigma$, we have $\sigma-eu=(1-e)\sigma$. Thus $(\sigma-eu)(y)=0$ for $y\in E$, while for $y\notin E$, the element $\sigma(y)$ is a nonunit of the local ring $f(y)$ and hence belongs to its maximal ideal. Since every maximal ideal of $A$ arises from a maximal ideal of a Pierce stalk,
we conclude that $\sigma-eu\in J$.

Thus every element of $A/J$ is the product of an idempotent and a unit. Hence $A/J$ is von Neumann regular and therefore has Krull dimension zero. Lemma~\ref{semiregular}\ref{semiregular 2} completes the proof.
\end{proof} 

Let $X\subseteq\Sigma(R)$ be a spectral space in the Zariski subspace topology, and let $X^{\mathrm{patch}}$ denote its intrinsic patch topology. Define the {\it $w$-localization} $X^{\rm wl}$ of $X$ to be  $$X^{\rm wl}=\{(S,T)\in X^{\mathrm{patch}}\times X:S\subseteq T\},$$ with the subspace topology, where the second factor has the Zariski topology. (The version we state here is homeomorphic to the   $w$-localization of $X$ constructed in \cite[Lemma~2.1.10]{BhattScholze}; see Lemma~2.5 and the proof of Lemma~7.13 in \cite{ScholzeDiamonds}.)  This space is $w$-local, and the  closed points of $X^{\rm wl}$ are the pairs $(S,S)$, which form a closed subspace  homeomorphic to $X^{\mathrm{patch}}$ \cite[Lemma 2.1.10]{BhattScholze}. The projection to $X$, $$ \pi\colon\; X^{\rm wl}\longrightarrow X\colon\; (S,T)\longmapsto T,$$ restricts on the set of closed points $C$  to a spectral map onto $X$. Under the canonical homeomorphism $X^{\mathrm{patch}}\cong C$, given by $S\mapsto(S,S)$, the restriction $\pi|_C$ corresponds to the identity map $X^{\mathrm{patch}}\to X$. Composing with the inclusion
$X\hookrightarrow\Sigma(R)$ gives the tautological patch bundle $\iota_X\colon\; X^{\mathrm{patch}}\to\Sigma(R)$. In the scheme case considered next, this identification allows us to compare global sections on \(X^{\mathrm{wl}}\) with the patch algebra of \(\iota_X\) by restricting sections to the closed-point subspace $C$.

Now suppose $(X,\O_X)$ is an integral separated quasi-compact scheme,  represented as in Example~\ref{scheme 1} so that $X \subseteq \Sigma(F)$, where $F$ is the function field of $X$.   
Then using the projection $X^{\rm wl} \rightarrow X$, we can pull back the 
 sheaf $\O_X$ to the sheaf $\pi^{-1}\O_X$ on $X^{\rm wl}$ to obtain a locally ringed space $(X^{\rm wl},\pi^{-1}\O_X)$.  As noted in the proof of the next proposition, it follows from results in \cite{BhattScholze} that this locally ringed space is an affine scheme. We show that this affine scheme is the prime spectrum of the associated patch algebra. 

\begin{theorem}\label{thm: scheme w-localization}
Let $X$ be an integral separated quasi-compact scheme, identified with a subspace of
$\Sigma(F)$ as in Example~\ref{scheme 1}, where $F$ is its function field.
Let $\pi\colon X^{\rm wl}\to X$ be its $w$-localization, and let $A=A(\iota_X)$ be the patch algebra of the tautological patch bundle associated with $X$.
Then there is an isomorphism of schemes, $$\Spec(A) \cong (X^{\rm wl},\pi^{-1}\O_X).$$  
\end{theorem}

\begin{proof}
By \cite[Lemma~7.14]{ScholzeDiamonds}, the spaces $U^{\rm wl}$, where $U$ ranges over an affine open cover of $X$, form an open cover of $X^{\rm wl}$. Using the  restrictions of $\pi^{-1}\mathcal O_X$ to the open sets in this cover, these are affine schemes by \cite[Lemma~2.2.4]{BhattScholze}, so $(X^{\rm wl},\pi^{-1}\mathcal O_X)$ is a scheme. Since its underlying space is $w$-local,
every open cover admits a refinement by pairwise disjoint
clopen subsets \cite[Definition~2.1.1 and Lemma~2.1.4]{BhattScholze}. Hence it is affine by \cite[Theorem~4.1]{Schroer}.

Let $Z=X^{\rm wl}$, and identify its closed-point subspace with $X^{\mathrm{patch}}$ as above. The first projection induces a continuous retraction $c\colon Z\to X^{\mathrm{patch}}\colon \;(S,T)\longmapsto S$,
where $c(z)$ corresponds to the unique closed specialization of $z$.
Since $X$ is integral, $\mathcal O_X$ is a subsheaf of the constant sheaf $\underline F_X$. Hence $\pi^{-1}\mathcal O_X\subseteq\underline F_Z$, with stalk $\mathcal O_{X,\pi(z)}$ at $z\in Z$ \cite[Tag~008H]{stacks-project}. 
The global sections of $\underline F_Z$ are locally constant $F$-valued functions. Since membership in the subsheaf $\pi^{-1}\mathcal O_X$ can be checked on stalks, its global sections are precisely the locally constant functions $\tau\colon\; Z\to F$ satisfying
$\tau(z)\in\mathcal O_{X,\pi(z)}$ for every $z\in Z$. Since each such function $\tau$ is locally constant, its fibers are closed, so the fiber containing $z$ contains $c(z)$ also, since $c(z)$ is the unique closed point in the closure of $\{z\}$. Since $\tau(z) = \tau(c(z))$, this shows that $\tau$ is uniquely determined by its values on the closed points of $Z$, and so two such functions $\tau$ that agree on closed points agree everywhere. 

Conversely, represent $A$ as in Theorem~\ref{C rep}. If $\sigma\in A$, then $\sigma\in C(X^{\mathrm{patch}},F_{\mathrm{disc}})$ and $\sigma(x)\in\mathcal O_{X,x}$ for every $x\in X$. The function $\sigma\circ c$ is locally constant on $Z$, and
\[
\sigma(c(z))\in\mathcal O_{X,c(z)}
\subseteq\mathcal O_{X,\pi(z)}
\]
for every $z\in Z$.
Thus $\sigma\circ c$ is a global section of $\pi^{-1}\mathcal O_X$.

This shows that the map that sends each global section of $\pi^{-1}\O_X$ to its restriction to $X^{\rm patch}$ induces an isomorphism $\Gamma(Z,\pi^{-1}\mathcal O_X) \cong A.$ Since $(Z,\pi^{-1}\mathcal O_X)$ is affine, this yields the claimed isomorphism of schemes.
\end{proof}

In \cite{Wiegand1971}, Wiegand shows that for every  ring $R$, there is a homomorphism $\eta\colon\; R\to T(R)$ to a von Neumann regular (i.e., absolutely flat) ring $T(R)$ such that every homomorphism from $R$ to a commutative von Neumann regular ring factors uniquely through $\eta$, and this property determines $T(R)$ up to a unique $R$-algebra isomorphism.

\begin{corollary} Let $D$ be an integral domain, and let $A$ be the patch algebra of the patch bundle $\Spec(D)^{\rm patch} \rightarrow \{D_P:P \in \Spec(D)\}$.  Then $A/J(A)$ is the universal von Neumann regular ring of $D$.
\end{corollary}

\begin{proof} Applying the preceding theorem to $X=\operatorname{Spec}(D)$ and using \cite[Lemma~2.2.4]{BhattScholze}, we obtain a canonical $D$-algebra isomorphism from $A$ to the ring of global sections of $(X^{\rm wl},\pi^{-1}\O_X)$. By \cite[Example~2.2.5]{BhattScholze}, the natural map $D\rightarrow A/J(A)$ is universal among homomorphisms from $D$ to  von Neumann regular rings, so  the corollary follows.
\end{proof}

This identification also highlights a connection between two sheaf-theoretic constructions. While the patch algebra $A$ arises from the local rings $D_P$ over $\Spec(D)^{\rm patch}$, Wiegand's universal regular ring is realized as the ring of global sections of a sheaf of their residue fields $\kappa(P)$ over the same Stone space \cite[Section~1]{Wiegand1971}. Thus, passing from $A$ to $A/J(A)$ connects these two constructions.

\section*{Acknowledgments}

Laura Cossu gratefully acknowledges the hospitality of the Department of Mathematical Sciences at New Mexico State University during her visits in Spring 2025 and Spring 2026, when much of this work took shape.
Bruce Olberding was supported by NSF grant DMS-2231414.

The authors thank Roger Wiegand for inspiring them to explore connections between patch algebras and his work on universal von Neumann regular rings.
\section*{AI disclosure}

This work was initiated in 2024. Its main mathematical results were developed collaboratively by the authors during Laura Cossu’s visits to New Mexico State University in Spring 2025 and Spring 2026, and in the course of regular online meetings. ChatGPT (GPT-6 Astra, OpenAI) was used only in the final stages of manuscript preparation, to improve the exposition and the organization of the paper, and to develop some illustrative examples. The authors take full responsibility for the mathematical content, the accuracy of the references, and the final version of the manuscript.

\end{document}